\documentclass[twoside, 11pt]{article}

\usepackage{amssymb, amsmath, latexsym, mathrsfs, verbatim, calc, graphicx}
\usepackage{colortbl}

\newtheorem{assumption}{Assumption}
\def\qed{ \ \vrule width.2cm height.2cm depth0cm\smallskip}

\newcommand{\ba}{\begin{array}}
\newcommand{\ea}{\end{array}}
\newcommand{\be}{\begin{equation}}
\newcommand{\ee}{\end{equation}}
\newcommand{\bea}{\begin{eqnarray}}
\newcommand{\eea}{\end{eqnarray}}
\newcommand{\beaa}{\begin{eqnarray*}}
\newcommand{\eeaa}{\end{eqnarray*}}

\def\dbR{{\mathop{\rm I\negthinspace R}}}
\def\a{\alpha}
\def\b{\beta}
\def\g{\gamma}
\def\d{\delta}
\def\e{\varepsilon}

\def\l{\lambda}

\def\si{\sigma}
\def\t{\tau}
\def\f{\varphi}
\def\th{\theta}
\def\o{\omega}
\def\h{\widehat}
\def\G{\Gamma}
\def\D{\Delta}

\def\L{\Lambda}

\def\cF{{\cal F}}

\def\cO{{\cal O}}

\def\cS{{\cal S}}

\def\cV{{\cal V}}

\def\cZ{{\cal Z}}
\def\no{\noindent}

\def\ss{\smallskip}
\def\ms{\medskip}
\def\bs{\bigskip}
\def\q{\quad}
\def\qq{\qquad}

\def\pa{\partial}

\def\cds{\cdots}

\def\bF{{\bf F}}

\def\qed{ \hfill \vrule width.25cm height.25cm depth0cm\smallskip}

\newcommand{\basa}{\begin{assumption}}
\newcommand{\easa}{\end{assumption}}

\newcommand{\bas}{\begin{assum}}
\newcommand{\eas}{\end{assum}}

\def\limsup{\mathop{\overline{\rm lim}}}

\def\pa{\partial}
\def\h{\widehat}

\def\cds{\cdots}

\def\sgn{\hbox{\rm sgn$\,$}}

\def\dis{\displaystyle}

\def\bF{{\bf F}}
\def\cad{c\`{a}dl\`{a}g~}

\def\one{{\hbox{\bf 1}}}

\newtheorem{thm}{Theorem}[section]
\newtheorem{lem}[thm]{Lemma}

\newtheorem{prop}[thm]{Proposition}
\newtheorem{rem}[thm]{Remark}
\newtheorem{eg}[thm]{Example}
\newtheorem{defn}[thm]{Definition}
\newtheorem{assum}[thm]{Assumption}

\begin{document}

\title{\bf Optimal Portfolio Selection under Concave Price Impact}
\author{Jin Ma\thanks{\noindent Department of Mathematics, University of
Southern California, Los Angels, CA 90089; email: jinma@usc.edu.
This author is supported in part by NSF grants \#DMS
0806017 and 1106853}., ~Qingshuo Song\thanks{ \noindent
Department of Mathematics, City University of Hong Kong; email:
songe.qingshuo@cityu.edu.hk. The research of this author is
supported in part by the Research Grants Council of Hong Kong No.
9041545.}, ~Jing Xu\thanks{ \noindent School of Economics and
Business Administration, Chongqing University, Chongqing 400030,
China; email: xujing8023@yahoo.com.cn.}, ~and ~Jianfeng
Zhang\thanks{ \noindent Department of Mathematics, University of
Southern California, Los Angels, CA 90089; email:
jianfenz@usc.edu. This author is supported in part by NSF grant
\#DMS 1008873.}
} 
\date{}
\maketitle

\begin{abstract}
In this paper we study an optimal portfolio selection problem under instantaneous price impact. Based on some empirical analysis in the literature, we
model such impact as a concave function of the trading size when the trading size is small. The price impact can be thought of as either a
liquidity cost
or a transaction cost, but the concavity nature of the cost leads to some fundamental difference from those in the existing literature. We show that the problem can be
reduced to an impulse control problem, but without fixed cost, and that the value function is a viscosity solution to a special type of Quasi-Variational Inequality
(QVI). We also prove directly (without using the solution to the QVI) that the optimal strategy exists and more importantly, despite the absence of a fixed cost,  it is still in a ``piecewise constant"
form, reflecting a more practical perspective.
\end{abstract}

\vspace{5mm}

\vfill



\eject

\section{Introduction}
\setcounter{equation}{0}

Modeling of the liquidity risk has attracted strong attention in the recent years in the quantitative
finance literature, and there have been numerous publications on the subject.  Among others, one of  the core issues is to understand the price impact of individual tradings.  Motivated by empirical observations, Bouchaud, Farmer, and Lillo
\cite{BFL}  (and the references therein) suggested a price impact model in which the trading size affects the price
in a ``concave" way, when the trading size is small. Such a (concavity) assumption apparently leads to some fundamental differences from many existing results (see more detailed discussion in \S2), and this paper is an attempt to understand these differences in the context of an optimal portfolio selection problem. Roughly speaking, we shall argue that under such  a
concavity assumption, the optimization problem can be reduced to an impulse control problem without a fixed cost, but nevertheless the optimal strategy still exists and, somewhat surprisingly, it is in a piecewise constant form. One can then
easily conclude that  the liquidity cost does exist.

Our model is mainly motivated by the work of Cetin, Jarrow, and Protter \cite{CJP},  in which the liquidity cost
was characterized by the so-called ``supply curve". The main feature of the model (along with its subsequent work by Cetin, Jarrow, Protter, and  Warachka \cite{CJPW}) is that the dependence of the supply curve on the trading size is essentially quadratic when the size is small. Furthermore, it is shown in \cite{CJP} that the supply-curve-based liquidity cost could be
eliminated if one is allowed to split any
(large) order into many small pieces, and trade them infinitely frequently (this amounts to saying that the continuous trading
is allowed).
Such a point was later amplified by Bank and Baum \cite{BB}, in which they
proved that  one can always approximate a trading strategy by
those that have continuous and finite variation paths, consequently the liquidity cost could always be eliminated.
But on the other hand, both empirical evidences and other theoretical studies indicate that the liquidity risk does exist,
even in the continuous trading paradigm. For instance, by considering the Gamma constraint on the admissible (continuous!)
portfolios
and by using the so-called second order backward SDEs, Cetin, Soner, and Touzi \cite{CST} proved that the super-hedging price is in general higher than the Black-Scholes price,  and thus  the liquidity
cost must exist.
Also, to make the model more realistic, various constraints on the trading strategies have been added in order to avoid the vanishing liquidity cost. For example, Cetin and Rogers \cite{CR} assumed that any two consecutive transactions have to be one unit of time apart.  In a different work, Ly Vath, Mnif, and Pham \cite{LMP}  assumed heavy liquidity cost if two transactions were made too closely. We should note, however, in the last two works
the optimal strategy being piecewise constant is (essentially) assumed  {\it exogenously}. 
The main message of
our result is that the concavity assumption of the liquidity cost provides an {\it endogenous} structure, from which the
optimal strategy becomes intrinsically ``piecewise constant", even in the absence of a fixed cost. 

It is worth noting that since all the liquidity costs mentioned above have instantaneous (or temporary) price impact, technically
they are equivalent to a type of transaction costs. Consequently, our approach can be easily applied to problems with transaction costs, which has been studied extensively (see, e.g., \cite{EH}, \cite{SS}, \cite{Korn}, \cite{OS},
\cite{Chang1}, \cite{Chang2}, \cite{LMP}, and the references therein). Most results in the literature  assume either fixed cost, or proportional cost,  or the linear combination of them. To be more precise, if we denote $c(z)$ to be the price impact or the transaction cost when the trading size is $z$, and we assume $c(z) \sim |z|^\a$ when $z$ is small, then the  fixed cost case corresponds to $\a = 0$, proportional cost or linear price impact case corresponds to $\a=1$, and the price impact in \cite{CJP} corresponds to $\a=2$. When $\a>1$, the liquidity (or transaction) cost vanishes in approximate sense by allowing multiple instantaneous trading. When $\a=1$, this is typically a singular control problem and the optimal strategy is continuous. When $\a=0$, this is typically an impulse  control problem and the optimal strategy is discrete. We essentially assume $0< \a <1$, which is consistent with the concavity of the price impact as observed in \cite{BFL}.  We show that our problem is essentially an impulse control problem, but possibly without fixed cost.

Our second goal in this paper is to prove the existence of the optimal strategy and argue that it must be piecewise
constant. We note that unlike most of impulse control problems in the literature, we do not assume that the
cost function is strictly positive (no fixed cost). Thus the HJB equation, being a quasilinear-variational inequality (QVI),
does not have  a smooth solution in general. Consequently, the construction of the optimal strategy, whence in many
cases the existence of it, become problematic if one follows the standard verification theorem approach
(cf., e.g. \cite{OS}).
In this paper we shall attack the existence of optimal strategy
directly. We first consider a sequence of approximating problems for
which the strategies are restricted to a fixed number (say, $n$)
of trades. We show that for each $n$ the optimal strategy, denoted by
$Z^n$, exists.
The main technical part in this analysis turns out to be some
uniform estimates on the number of jumps of $Z^n$. These estimates
will enable us to study the regularity of the value function and to
construct the optimal strategy. We should note that the regularity
of the value function, which we need to construct the optimal
strategy, is weaker than those that are commonly seen in the
literature.

The rest of the paper is organized as follows. In Section 2 we
formulate the problem and state the main result. In Sections 3 and 4 we study the
approximating value function $V^n$ and its corresponding optimal
strategy $Z^n$. In Section 5 we obtain uniform estimates of $Z^n$,
which leads to the regularity of the value function $V$.
In Section 6  we study the optimal strategy of the
original problem.
Finally in Section 7 we give some technical proofs.

\section{Problem Formulation}
 \setcounter{equation}{0}

\subsection{The model}
Let $(\Omega, \cF, P;\bF)$ be a complete filtered probability space
on a finite time interval $[0,T]$ and $W$ be a standard Brownian motion.
We assume that the filtration $\bF=\{\cF_t\}_{t\ge 0}$ is generated by $W$, augmented by all the $P$-null sets as usual.
The financial market consists of  two assets, a bank account and a stock. For simplicity, we
assume that the interest rate is $0$. Let $X$ denote the fundamental value of the stock which follows the stochastic differential equation:
 \bea
 \label{X}
d X_t = b(t,X_t)dt + \sigma(t,X_t)dW_t.
 \eea
In this paper we consider the liquidity cost in the following general form:  if one buys $z$ shares of the stock (sells $-z$ shares if $z<0$) at  time $t$, then the liquidity cost of the trade is $c(t,X_t,z)$, where $c$ is a deterministic function satisfying 
$c(t,x,0)=0$; and 
\bea
\label{cproperty}
\mbox{$c$ is increasing in $z$ when $z>0$ and decreasing in $z$ when $z<0$.}
\eea
We shall give more specific assumptions on the cost function $c$  in the next subsection. But we remark here that
if $c_0 := \inf_{(t,x), z\neq 0} c(t,x, z)>0$, then $c_0$ represents a ``fixed cost". The following example shows that such
a positive lower bound usually does not exist in the context of liquidity cost.

\begin{eg}{\rm
\label{ExCJP} Consider the ``supply curve" $\cS(t,X_t,z)$ defined in  \cite{CJP}, in which  $X_t$ is the fundamental price and  $z$ is the trading size at time $t$. We can view $\cS$ as the market price of the stock, satisfying 
\bea
\label{cSproperty}
\cS(t, X_t, 0) = X_t, ~\mbox{and $\cS$ is increasing in $z$.}
\eea
Thus the liquidity cost should naturally be defined by
\bea
\label{cSc}
c(t,X_t,z) := z[\cS(t,X_t,z) - X_t].
\eea
One can easily check that the $c$ satisfies (\ref{cproperty}).
\qed}
\end{eg}
We remark that in Example \ref{ExCJP}, if  $\cS$ is smooth in $z$, then  $c(t,X_t,z)\sim z^2$ when $z$ is small. Namely $z\mapsto c(t,X_t,z)$ is {\it convex} for $z$ small.
In this paper, however,  we are interested in the case where $c(t, X_t,z) \sim |z|^\a$ for some $0< \a<1$, as supported by \cite{BFL}. Therefore it is fundamentally different from the case  in \cite{CJP}.

We next  consider admissible trading strategies $Z$. We assume $Z$ is $\bF$-adapted, \cad, and piecewise constant.  Let $\tilde Y$ denote the total value invested in the riskless  asset, and define
$
Y:=\tilde Y + ZX.
$
Assuming that the interest rate is $0$, then except  for countably many $t\in D_Z:= \{t\in [0,T]: \d Z_t := Z_t - Z_{t-}\neq 0\}$,
one has
\bea
\label{SF1}
d \tilde Y_t = 0 &\mbox{and thus}& d Y_t = Z_t dX_t.
\eea
Namely, $Z$ is ``self-financing". Furthermore, for $t\in D_Z$ (i.e., $\d Z_t\neq 0$), we impose the standard self-financing constraint:
\bea
\label{SF2}
 \d Y_t + c(t, X_{t-}, \d Z_t)=\d Y_t+c(t,X_t,\d Z_t)=0.
\eea
We note that (\ref{SF2}) simply means that no instantaneous profit can be made by changing the investment positions. In the
case of supply-curve (Example \ref{ExCJP}), the equation (\ref{SF2}) amounts to saying that (noting that $X$ is continuous)
\beaa
\d \tilde Y_t  +\d Z_t X_t +c(t, X_t, \d Z_t)  = \d \tilde Y_t+\d Z_t \cS(t, X_t, \d Z_t) =0.
\eeaa
This is exactly the standard idea of ``budget constraint".

\subsection{The optimization problem}
We now introduce our optimization problem on a subinterval $[t,T]$. 
Let $X^{t,x}$ denote the solution to SDE (\ref{X}) with initial value $X_t=x$, a.s. Given $(x,y,z)$ and an admissible trading strategy $Z$,  we shall set $Y_{t-} := y$ and $Z_{t-}:=z$. Then by (\ref{SF1}) and (\ref{SF2}) we have
\bea
\label{Y}
Y^{t,x,y,z,Z}_T := Y_T = y + \int_t^T Z_s dX_s^{t,x} - \sum_{t\le s\le T} c(s, X_s^{t,x}, \d Z_s).
\eea
Let $U$ be a  terminal payoff function, then our optimization problem is:
\bea
\label{V}
V(t,x,y,z) := \sup_{Z\in \cZ_t} E\Big[ U(Y^{t,x,y,z,Z}_T)\Big].
\eea
Here the set $\cZ_t$  of the admissible strategies is defined rigorously at below:

 \begin{defn}
\label{admissible} Given $t\in [0,T]$, the set of
admissible strategies, denoted by $\cZ_t$, is the space of
$\bF$-adapted processes $Z$ over $[t,T]$ such that, for a.s. $\o$,

(i) $Z$ is  \cad and piecewise constant with finitely many jumps;

(ii) $Z_T=0$, and $|Z|\le M$.
\end{defn}
It is worth noting that a key assumption in Definition \ref{admissible} is that $Z$ is piecewise constant and has only finitely many jumps. While this is obviously
more desirable in practice, it is by no means clear that an optimal strategy can be found in such a form. Thus the main
goal of this paper is to show that the concavity assumption on $c$, see  {\bf (H4)} below,  implies the existence of an optimal strategy in $\cZ_t$.

\begin{rem}
{\rm
(i) We require $Z$ to be \cad for notational convenience. One can easily change it to c\`agl\`ad if necessary.

(ii) Due to the liquidity risk, if $Z_T\neq 0$, the payoff of $Y_T$ is not clear.  As in \cite{CJP} and \cite{CST}, we require $Z_T=0$ so that $Y_T=\tilde Y_T$. An alternative way is to introduce a payoff function $U(\tilde Y_T, Z_T)$ on both accounts, see, e.g. \cite{CPT} in the formulation of superhedging.

(iii) The assumption that $Z$ is bounded is merely technical. This
restriction can be removed, with some extra efforts on the estimates,  by requiring that the cost function
$c$ satisfies certain growth condition, for example, $\lim_{|z|\to\infty}\inf_x |c(z,x)|\slash |z|=\infty.$
We prefer not to pursue such complexity in this paper. In fact, we will impost some stronger technical assumptions
in order not to distract our attention from the main focus of the paper.
\qed
}\end{rem}

\begin{rem}
\label{rem-Z}
{\rm
Technically, the optimization problem (\ref{V}) can be extended to the cases where admissible strategies
are allowed to be general $\bF$-adapted, \cad
processes. But in that case we need to redefine the aggregate
liquidity cost. For example, we can consider the aggregate
cost in the following forms:
 \bea
\label{agtrans}
 \sup_\pi \sum_{i=0}^\infty
 c(\t_i, X_{\t_{i}}, Z_{\tau_i} - Z_{\t_{i-1}}), \q \mbox{\rm or}\q \lim_{|\pi|\to 0}
 \sum_{i=0}^\infty c(\t_i, X_{\t_{i}}, Z_{\tau_i} - Z_{\t_{i-1}}),
 \eea
where the supreme is over all possible random partitions of $[t,T]$
$\pi$: $t=\t_0<\t_1<\cds \le T$; 
and $|\pi|$ is the ``mesh size" of
the partition. Then, under our conditions in next subsection on the function
$c$ , one can show that the value function $V$ would be the same as the one
where the supreme is taken over only piecewise constant strategies.
Namely, it suffices to consider only $\cZ_t$, and thus the
aggregate cost (\ref{agtrans}) is again reduced to that in
(\ref{Y}).

 However, for more general $c$, typically there is no optimal strategy in $\cZ_t$ and then one has to extend the space to allow more complex strategies. The following two cases are worth noting.

(i) Assume that $c(t,x,z)=|z|$. Then $\sup_\pi \sum_{i=0}^\infty
 c(\t_i, X_{\t_{i}}, Z_{\tau_i} - Z_{\t_{i-1}})=\int_t^T|dZ_r|$, the total
 variation of the process $Z$. This problem then becomes a more or less
 standard singular (or impulse) stochastic control problem (cf.
 e.g., \cite{EH}, \cite{Korn}, and \cite{LMP}). In these cases the
 optimal controls are of bounded variation, but not necessarily piecewise
 constant.

 (ii) Assume  the supply curve $\cS(t,x,z)$ is smooth,  as proposed in
\cite{CJP} and \cite{CJPW}. Then $c(t,x,z) \sim  z^2$ when $z$ is small. For any (random) partition $\pi:
 t=\t_0<\t_1<\cds\le T$ and any $\bF$-adapted semimartingale $Z$ satisfying $Z_T=0$, we have
 \beaa
&& \sum_{i=0}^\infty  c(\t_i, X_{\t_{i}}, Z_{\tau_i} - Z_{\t_{i-1}})= \sum_{i=0}^\infty  [Z_{\tau_i} - Z_{\t_{i-1}}][\cS(\t_i, X_{\t_{i}}, Z_{\tau_i} - Z_{\t_{i-1}}) - X_{\t_i}]\\
&=&\sum_{i=0}^\infty  [Z_{\tau_i} - Z_{\t_{i-1}}][\cS(\t_i, X_{\t_{i}}, Z_{\tau_i} - Z_{\t_{i-1}}) - \cS(\t_i, X_{\t_{i}}, 0)]\\
 &\to& \sum_{t\le s\le T} \d Z_s[\cS(s, X_s,\d Z_s)-\cS(s,X_s, 0)]+\int_t^T \frac{\pa
 \cS}{\pa z}(s, X_s,0)d[Z,Z]^c_s.
 \eeaa
 This recovers the liquidity cost in \cite{CJP} and \cite{CJPW}, and
 in this case it is natural to set the admissible strategies as semimartingales.
\qed
}
\end{rem}

\subsection{Technical assumptions}

We now present our technical conditions. As mentioned in Remark \ref{rem-Z}, our main focus is to show 
that the concavity assumption on $c$ implies the existence of an optimal strategy in $\cZ_t$. However, in order not to 
over complicate our estimates, we shall impose some stronger technical conditions, some of which 
might be more than necessary. We remark that our approach can be extended to more general cases.

We first assume that all processes in this paper are
one dimensional and, as mentioned already, the interest rate is $0$. Moreover,  we shall make use of the following assumptions:
\begin{itemize}

 \item[{\bf (H1)}] The coefficients $b$ and $\si$ in (\ref{X}) are bounded and
 uniformly Lipschitz continuous in $x$, with a common Lipschitz
 constant $K>0$.

\item[{\bf (H2)}] The terminal payoff function $U$ is concave,
increasing such that
 $0<\l\le U'\le \L$ on $(-\infty, \infty)$ for some constants $0<\l<\L$.

 \item[{\bf (H3)}] The cost function $c$ depends only on the trading size $z$, and satisfies:

 (i)  $c(0)=0$ and $c(z)> 0$ for all $z\neq 0$;.

(ii) $c$ is increasing  in $[-2M,0)$ and decreasing in $(0,2M]$; moreover, in both intervals,  $c$ is uniformly continuous with
 the same modulus of continuity function $\rho$.

 (iii) the following subadditive property holds:
 \bea
 \label{subadd}
 c(z_1+z_2)\le c(z_1)+c(z_2), ~\mbox{for any $z_1,z_2$ such that $|z_1|, |z_2|,
    |z_1+z_2|\le 2M$}.
\eea
%
 \item[{\bf (H4)}] There exists a constant $\e_0>0$   such that

 (i)  $c$ is  concave in $(0, 2\e_0]$ and in $[-2\e_0, 0)$, and
 \bea
 \label{c0infty}
\eta_\th :=  \limsup_{z\to 0}{ c(\th z) \over  c(z)}  < \th, ~ \mbox{for}~ \th = {3\over 2}, 2, 3, ~~\mbox{and}~~ \g :=  \limsup_{z\to 0}{ c(-2z) -c(-z)\over  c(z)}
<\infty.
  \eea

  (ii) $c$ is uniformly Lipschitz continuous in $[-2M, -{\e_0}]\cup [{\e_0}, 2M]$ with a Lipschitz constant $L_0$.
  \end{itemize}

We conclude this subsection by several important remarks.
\begin{rem}
\label{rem-Inada}
{\rm The assumption ({\bf H2}) indicates that the terminal payoff $U$ is essentially a ``utility
function", except that it violates the well-known {\it Inada} condition:
\bea
\label{Inada}
\lim_{y\to -\infty} U'(y)= \infty, \q \lim_{y\to\infty} U'(y) = 0.
\eea
This is mainly for technical simplifications. The following observations are worth noting.

(i) If there is a fixed cost, namely if the cost function $c$ satisfies
\bea
\label{fixcost}
 c(z)\ge c_0>0 ~~\mbox{for all}~~ z\neq 0,
 \eea
 then one can prove our main result Theorem \ref{thm-main} under Inada condition 
 (\ref{Inada}) (see also Remark \ref{rem-c}-(iii) below). In fact, in this case the conditions on $c$ can also be further relaxed.

(ii) In the case when $U(y) = - e^{-y}$, $c(z) = |z|^\a$ for some $0<\a<1$, and $b(t,x) = b_0$, $\si(t,x) = \si_0$, then the
assumptions ({\bf H1}), ({\bf H3}), ({\bf H4}), and (\ref{Inada}) are all satisfied, one can easily check that $V(t,x,y,z) = -e^{-y} \cV(t,z)$, where  
\bea
\label{calV}
\cV(t,z) := \inf_{Z\in \cZ_t} E\Big[\exp\big(- b_0\int_t^T Z_s ds - \si_0 \int_t^T Z_sdW_s + \sum_{t\le s\le T} |\d Z_s|^\a\big)\Big]. 
\eea
Thus the optimization problems (\ref{V}) and (\ref{calV}) are equivalent.
By utilizing the structure of $\cV$ and modifying our arguments slightly we can also prove our main result in this case.

We believe our results hold true  under even more general conditions. However, since the main focus of this paper is the impact of the concave cost function $c$, we choose not to over-complicate this already lengthy paper, and content ourselves
with the stronger condition ({\bf H2}) instead.
\qed}
\end{rem}

\begin{rem}
\label{rem-c}  {\rm (i) We require the concavity of $c$ only around $0$. Typically, $c$ is convex when $z$ is large, as in the standard literature of liquidity risk.

(ii)  The typical case satisfying {\bf (H3)} and {\bf (H4)} is:  $c(z)=c_0|z|^{\a}$, $0<\a<1$. The condition (\ref{c0infty}) captures the behavior of $c$ around $0$. We consider those three values of $\th$ just for technical reasons. One can of course make the assumption more symmetric by strengthening the condition to $\eta(\th) <1$ for all $\th>1$. The assumption on $\g$ is merely technical. However, one cannot remove (\ref{c0infty}) for free. For example, $c(z) = |z|$ violates (\ref{c0infty}) and we know in this case the optimization problem becomes a singular control problem, see Remark \ref{rem-Z} (i).

(iii) Another typical case is when there is a fixed cost, namely (\ref{fixcost}) holds.
Since in this case (\ref{c0infty}) automatically holds, we do not need the concavity assumption in {\bf (H4)} and our main results will still be valid. See Theorem \ref{thm-main} below.

 (iv) Note that we allow $c(0+)>0$ and/or $c(0-) >0$ in (H4).  Moreover, combing the arguments for the two cases in (i) and (ii), we can easily prove our results in the case that (H4) holds in $(0, 2\e_0]$ and $c(z) \ge c_0 >0$ for $z<0$, and the case that   (H4) holds in $[-2\e_0, 0)$ and $c(z) \ge c_0 >0$ for $z>0$.
\qed}
 \end{rem}

\begin{rem}
\label{rem-subadd}
{\rm (i) In this remark we justify the subadditive property  (\ref{subadd}). Note that our goal is to solve (\ref{V}). For general $c$, by possibly splitting a transaction into many small pieces, we define,
 $$
 \tilde c(z) := \inf\{c(z_1)+\cds+c(z_n): |z_i|\le 2M, z_1+\cds+z_n=z,\forall n\}.
 $$
Then it is easy to see that $\tilde c\le c$ and $\tilde c$ satisfies
(\ref{subadd}). Replacing $c$ by $\tilde c$ in (\ref{Y}) we have
 $$
 \tilde Y_T := y + \int_t^T Z_s dX_s  - \sum_{t\le s\le T}\tilde c(\d Z_s);\q \tilde V(t,x,y,z):= \sup_{Z\in\cZ_t} E \Big[U(\tilde Y_T) \Big].
 $$
Under the continuity of $U$, one can easily show that $\tilde V = V$.
In other words, we can always replace the cost function
$c$ to one that satisfies (\ref{subadd}).

(ii) If the cost function $c$ satisfies
$c(z)\le C|z|^\a$ for some constants $C>0$ and $\a>1$ near $z=0$,
then the corresponding $\tilde c(z)\equiv 0$. To see this, note that
for any $z$ and large $n$ we have
 $$
 \tilde c(z) \le \sum_{i=1}^n c({z\over n}) \le C\sum_{i=1}^n |{z\over n}|^\a
 \le {CM^\a\over n^{\a-1}}\to 0,\q \mbox{\rm as~} n\to\infty.
 $$
Thus the optimization problem is reduced to a standard one without
liquidity cost. This is consistent with the result of \cite{CJP}, where $\a=2$.
\qed}
\end{rem}

\subsection{Main result}
For any $Z\in \cZ_t$,  we shall always denote
\bea
\label{ti}
\t_0:= t,\q \t_i:= \inf\{s>\t_{i-1}: Z_s \neq
Z_{\t_{i-1}}\}\wedge T,~~ i=1,\cds
\eea
Then clearly $\t_i < \t_{i+1}$ whenever $\t_i<T$, $\t_i = T$ when $i$ is large enough, and
\bea
 \label{Zcont}
 Z_s = \sum_{i=1}^\infty Z_{\t_{i-1}}\one_{[\t_{i-1},\t_i)}(s), \qq
 s\in[t,T].
 \eea
Recall that $Z_{t-}=z$. Let $N(Z)$ denote the number of jumps of $Z$, that is,
\bea
\label{NZ}
N(Z) := \sum_{t\le s\le T} {\bf 1}_{\{\d Z_s\neq 0\}} = \sum_{i=0}^\infty {\bf 1}_{\{Z_{\t_i} \neq Z_{\t_{i-1}}\}}.
\eea
Our main result of the paper is:

 \begin{thm}
 \label{thm-main}
 Assume {\bf (H1) -- (H3)}, and assume either (\ref{fixcost})  or {\bf (H4)}
 is in force. Then for any $(t,x,y,z)$,  the optimization problem (\ref{V}) admits an optimal strategy $Z^* \in \cZ_t$. Moreover, $E[N(Z^*)] <\infty$.
 \end{thm}

 \section{The Approximating Problems}
  \setcounter{equation}{0}

In this section, we shall approximate the
original optimization problem (\ref{Y}) and (\ref{V}) by those with
only fixed number of transactions, for which the optimal strategies
are easier to find.
To begin with, for any $n\ge 1$ we consider a reduced problem with at most $n$ transactions:
 \bea
 \label{Vn}
 V^n(t,x,y,z) := \sup_{Z\in\cZ_t^n(z)} E\{U(Y^{t,x,y,z,Z}_T)\}~~\mbox{where}~~\cZ_t^n(z):= \{Z\in \cZ_t:
N(Z)\le n\}.
 \eea
 We note that, for $Z\in\cZ^n_t(z)$, if $Z_t = z$, then $\t_n=T$, and if $Z_t \neq z$, then $\t_{n-1}=T$. Moreover, when $n=1$, we have $\cZ^1_t(z) = \{z\one_{[t,\t)}\}$ for all stopping time $\t$, and
 \bea
 \label{V1}
 V^1(t,x,y,z) = \sup_{\t\ge t}E\Big\{U\big(y+z(X_\t^{t,x}-x)-c(-z)\big)\Big\}.
 \eea
It is then readily seen,  assuming (H1)--(H3),  that
 \bea
 \label{V1bound}
 |V^1(t,x,y,z)|\le C[1+|y|], \q  (t,x,y,z)\in[0,T]\times \dbR^2\times [-M,M].
 \eea
Here and in the sequel $C>0$ is a generic constant depending only on
$T, M,\l,\L, K$, and $|U(0)|$ in (H1)--(H3),
as well as $\sup_{|z|\le 2M}c(z)$, and it is allowed to vary from
line to line. 
 \begin{prop}
 \label{Vnconv}
 Assume (H1)--(H3). Then $V^n(t,x,y,z)\uparrow V(t,x,y,z)$, as $n\to\infty$; and
 \bea
 \label{vngrow}
  V^n(t,x,y,z)
  \le V(t,x,y,z) \le C[1+|y|], \q (t,x,y,z)\in[0,T]\times
 \dbR^2\times [-M,M].
 \eea
 \end{prop}
 {\it Proof.} It is clear by definition that  $V^n$ is increasing and $V^n\le V$.
We first show that (\ref{vngrow}) holds for $V$ (whence for $V^n$ as
well). For any $Z\in \cZ_t$, let us denote $X=X^{t,x}$ and
$Y=Y^{t,x,y,z,Z}$ for simplicity. Since  the liquidity cost is positive, we
have
 $$
 Y_T \le y + \int_t^T Z_s dX_s = y+\int_t^T Z_s b(s,X_s)ds +
 \int_t^T Z_s \si(s,X_s)dW_s.
 $$
Then, using the monotonicity of $U$ and boundedness of $b$, $\si$ and $Z$, we have
 \bea
 \label{EUYT}
 EU(Y_T)\le E\Big\{U(y+\int_t^T Z_s dX_s)\Big\}\le |U(0)|+\L
 \Big\{|y|+E\Big|\int_t^T Z_s dX_s\Big|\Big\}
 \le C[1+|y|].
 \eea
Since $Z$ is arbitrary, we prove (\ref{vngrow}).

We now show that $V^n\to V$, as $n\to\infty$. We first note that
$V^n$ is non-decreasing, and bounded from above, thanks to
(\ref{vngrow}). Thus $V^\infty(t,x,y,z) :=
\lim_{n\to\infty}V^n(t,x,y,z)$ exists,  and $V^\infty(t,x,y,z)\le
V(t,x,y,z)$, for all $(t,x,y,z)$. We need only show that
$V^\infty\ge V$. To this end, for any $Z\in \cZ_t$ we define
$Z^n_s := Z_s 1_{\{s< \t_{n-1}\}}$, $s\in[t,T]$. 
Clearly,
$Z^n\in\cZ^{n}_t(z)$. Denote 
  $Y^n:= Y^{t,x,y,z,Z^n}$. Then by the subadditivity assumption (\ref{subadd}) we have
 \bea
 \label{dYn}
 Y_T - Y^n_T = \int_{\t_{n-1}}^T
 Z_s dX_s- \sum_{i\ge n} c(\d Z_{\t_i})
 + c(-Z_{\t_{n-1}})
 \le  \int_{\t_{n-1}}^T Z_s dX_s. 
 \eea
Now, for any $n$, using (H2), (\ref{vngrow}), and (\ref{dYn}) we
have
 \bea
 \label{EUYTn}
E\{U(Y_T)\} & =& E\{U(Y^n_T)\} + E\Big\{U(Y_T)-U(Y^n_T)\Big\}\nonumber\\
&=& E\{U(Y^n_T)\} + E\Big\{\Big[\int_0^1U'(Y^n_T+\th(Y_T-Y^n_T))d\th\Big][Y_T - Y^n_T ]\Big\}\\
&\le& V^\infty(t,x,y,z) +  \L E\Big\{\Big|\int_{\t_{n-1}}^T Z_s 
dX_s\Big| 
\Big\}. \nonumber
 \eea
Next, Definition \ref{admissible} (iii) implies that
 $
 \lim_{n\to\infty}\Big\{\Big|\int_{\t_n}^T Z_s dX_s^{t,x}\Big| 
 \Big\}=0,$ $P$-a.s.
This 
enables us to let $n\to\infty$ in
(\ref{EUYTn}) and apply the Dominated Convergence Theorem to get
$ E\{U(Y_T)\} \le  V^\infty(t,x,y,z)$. Since this is true for any
$Z\in\cZ_t$, we conclude that $V(t,x,y,z)\le V^\infty(t,x,y,z)$,
proving the proposition.
 \qed

The next  result concerns the {\it uniform regularity} of $\{V^n:
n\ge 1\}$.
  \begin{prop}
  \label{regxy}
  Assume (H1)--(H3). Then, for any $n$, it holds that
  \bea
  \label{regx}
 |V^n(t,x_1,y,z)-V^n(t,x_2,y,z)|&\le&  C|\D x|;\\
  \label{regy}
 \l \D y \le V^n(t,x,y_1,z)-V^n(t,x,y_2,z)&\le& \L\D y,\q \forall y_1\ge y_2;\\
  \label{regt}
 |V^n(t_1,x,y,z)-V^n(t_2,x,y,z)|&\le& C|\D t|^{1\over 2}.
  \eea
Here and in the sequel, $\D \xi:= \xi_1-\xi_2$, $\xi=t,x,y,z$,
respectively.

Moreover, for $z_1>z_2>0$ or $z_1<z_2<0$, we have
\bea
  \label{Vnregz}
  -C[|\D z|+\rho(|\D z|)]\le V^n(t,x,y,z_1)-V^n(t,x,y,z_2)\le C[|\D z|+\rho_n(|\D z|)],
  \eea
  where $\rho$ is the modulus of continuity of $c$ in (H3) (iii), and
 \bea
 \label{rhon}
  \rho_n(|\D z|) := \sup\Big\{ \sum_{i=1}^n \rho(\th_i |\D z|): \th_1,\cds,\th_n\ge 0, \sum_{i=1}^n \th_i = 1\Big\}\le n\rho(|\D z|).
  \eea

 \end{prop}

In this below,  we present the proof of \eqref{regx}, \eqref{regy}, and
\eqref{regt} only.  The proof of (\ref{Vnregz}) is more involved and thus is  relegated to  Section 7.

  {\it Proof.}
First let us denote $X^i:= X^{t,x_i}$, $i=1,2$, and $\D X :=
X^1-X^2$. Then by the standard arguments in SDEs we know that
 \bea
 \label{dXest}
  E\{\sup_{s\in [t,T]}|\D X_s|^2\} \le C|\D x|^2.
  \eea
Next, for any $Z\in\cZ^n_t(z)$, denote $Y^i:= Y^{t,x_i, y,z,Z}$,
$i=1,2$, and $\D Y := Y^1-Y^2$. Then
  \beaa
  |\D Y_T|\le \int_t^T |Z_s||b(s,X^{t,x_1}_s) - b(s,X^{t,x_2}_s)|ds
  + \Big|\int_t^T Z_s [\si(s,X^{t,x_1}_s) - \si(s,X^{t,x_2}_s)]dW_s\Big|.
 \eeaa
Since $b$ and $\si$ are Lipschitz continuous and $Z$ is bounded,
(\ref{dXest}) leads to that
 \beaa
 \Big|E\Big\{U(Y^1_T) - U(Y^2_T)\Big\}\Big|^2
 \le C E\Big\{|\D Y_T|^2\Big\} \le C E\Big\{\int_t^T |Z_s\D X_s|^2ds\Big\}
 \le C|\D x|^2.
 \eeaa
Since $Z$ is arbitrary, (\ref{regx}) follows easily.

To prove (\ref{regy}) we denote, for any $Z\in\cZ^n_t(z)$ and
$y_1>y_2$, $Y^i:= Y^{t,x,y_i,z,Z}$, $i=1,2$, and $\D Y := Y^1-Y^2$. Note that $\D Y_T=\D y$, we have
$$
E\Big\{U(Y^1_T) - U(Y^2_T)\Big\} =
E\Big\{\Big[\int_0^1U'(Y^1_T+\th\D y)d\th\Big] \D y\Big\}.
$$
Thus (\ref{regy}) follows from (H2) immediately.

We next prove (\ref{regt}). Assume $t_1<t_2$. It is then
standard to show that
 \bea
 \label{Xt12}
 E\Big\{|X^{t_1,x}_t - X^{t_2,x}_t|^2\Big\}\le C|\D t|,\q  t\ge t_2>t_1.
 \eea
Now for any $Z\in\cZ^n_{t_2}(z)$, define $\tilde Z_t:=
z1_{[t_1,t_2)}(t) + Z_t 1_{[t_2,T]}$. Then $\tilde
Z\in\cZ^{n}_{t_1}(z)$. Denote $X^i:= X^{t_i,x}$, $i=1,2$, and
$Y^2:= Y^{t_2,x,y,z,Z}$, $\tilde Y^1=Y^{t_1,x,y,z,\tilde Z}$, then
\beaa
&& Y^2_T-\tilde Y^1_T = \int_{t_2}^T Z_t dX^2_t - \int_{t_1}^T \tilde Z_t
 dX^1_t\\
 &&= -z[X^1_{t_2}-x] + \int_{t_2}^T
 Z_t[b(t,X^2_t)-b(t,X^1_t)]dt +\int_{t_2}^T
 Z_t[\si(t,X^2_t)-\si(t,X^1_t)]dW_t
\eeaa
Now by standard arguments one can easily derive from (\ref{Xt12}) that
$$
EU(Y^2_T)- V^n(t_1,x,y,z) \le
 E\Big\{U(Y^2_T)-U(\tilde Y^1_T)\Big\}
 \le CE\Big\{|Y^2_T-\tilde Y^1_T|\Big\}\le C|\D t|^{1\over 2}.
 $$
 Since $Z\in\cZ^n_{t_2}(z)$ is arbitrary, we get
 \be
 \label{regt1}
   V^n(t_2,x,y,z)- V^n(t_1,x,y,z)\le C|\D t|^{1\over 2}.
   \ee

On the other hand, for any $\dis Z=\sum_{i=1}^n
Z_{\t_{i-1}}\one_{[\t_{i-1},\t_i)}\in \cZ^n_{t_1}(z)$, it is obvious that $Z\in \cZ^n_{t_2}(z)$.
Denote $Y_i := Y^{t_i,x,y,z,Z}$ and assume $\t_{j}\le t_2< \t_{j+1}$. Note that $Z_{t_2}= Z_{\t_j}$. Then,  by the subadditivity assumption (\ref{subadd}),
\beaa
Y^1_T - Y^2_T &=& \int_{t_1}^T Z_t dX_t^1-\int_{t_2}^T Z_t dX_t^2 - \sum_{i=0}^j c(\d Z_{\t_i}) + c(Z_{\t_j}-z)\\
&\le& \int_{t_1}^{t_2} Z_t dX_t^1+  \int_{t_2}^T Z_t dX_t^1 -\int_{t_2}^T Z_t dX_t^2\\
\eeaa
Since $b, \si$ and $Z$ are bounded, one can easily check that
\beaa
E\Big[\Big| \int_{t_1}^{t_2} Z_t dX_t^1\Big|\Big] =E\Big[\Big| \int_{t_1}^{t_2} Z_t [ b(t, X^1_t)d t + \si(t, X^1_t) dW_t]\Big|\Big] \le C|\D t|^{1\over 2}.
\eeaa
Moreover, note that $X^1_t = X^{t_2, X^1_{t_2}}_t$ for $t\ge t_2$.  Following the arguments for (\ref{regx})  we have
\beaa
&&E\Big[\Big|  \int_{t_2}^T Z_t dX_t^1 -\int_{t_2}^T Z_t dX_t^2\Big|\Big] \\
&=&E\Big[\Big| \int_{t_2}^T Z_t\Big[[b(t, X^1_t) - b(t, X^2_t)] dt + [\si(t, X^1_t) - \si(t, X^2_t)] dW_t\Big|\Big] \\
&\le&  CE\Big[|X^{t_1,x}_{t_2} - x|\Big] = CE\Big[\Big| \int_{t_1}^{t_2} [ b(t, X^1_t)d t + \si(t, X^1_t) dW_t]\Big|\Big] \\
& \le& C|\D t|^{1\over 2}.
\eeaa
Then, by the assumption (H2) on the payoff function $U$,
   \beaa
&&    E\Big\{U(Y^1_T)\Big\} - V^n(t_2,x,y,z)\le
    E\{U(Y^1_T)-U(Y^2_T)\}\\
     &\le& CE\Big\{\Big| \int_{t_1}^{t_2} Z_t dX_t^1\Big| + \Big|  \int_{t_2}^T Z_t dX_t^1 -\int_{t_2}^T Z_t dX_t^2\Big|\Big\}\le C|\D t|^{1\over 2}
 \eeaa
 Since $Z\in\cZ_{t_1}^n(z)$ is arbitrary, we get
 $
  V^n(t_1,x,y,z)- V^n(t_2,x,y,z)\le C|\D t|^{1\over 2},
  $
  which, together with (\ref{regt1}), implies (\ref{regt}).
  \qed

  We will also need the following result in next section.
 \begin{prop}
 \label{regz0}
  Assume (H1)-(H3). Then for any $n$ and any $(t,x,y)$,
  $$
  V^n(t,x,y,0+)\le V^n(t,x,y,0);\q V^n(t,x,y,0-)\le V^n(t,x,y,0).
  $$
  \end{prop}
  {\it Proof.} First by (\ref{Vnregz}) we know $V^n(t,x,y,0+)$ and $V^n(t,x,y,0-)$ exist.

  For $z>0$ and $Z^1=\sum_{i=1}^n Z^1_{\t_{i-1}}\one_{[\t_{i-1},\t_i)} \in\cZ_t^n(z)$, we define $Z^2  \in\cZ^n_t(0)$ as follows. Let $k:= \inf\{i: Z^1_{\t_i}\le 0\}$. We note that $k\le n$ since $Z_{\t_n}=0$. Define $Z^2_s := [Z^1_s - z]\vee 0$ for $s< \t_k$ and $Z^2_s:= Z^1_s$ for $s\ge \t_k$. Denote $\Delta Z := Z^1 - Z^2$.
It is  straightforward to check that
 $$
 0\leq \Delta Z_{\tau_i}\leq z ~~\mbox{and}~~ \d Z^1_{\tau_i} \d Z^2_{\tau_i}\ge 0,\q  i=0,\cdots n.
 $$
  Note that
   $$
   Y^{t,x,y,z,Z^1}_T-Y^{t,x,y,0,Z^2}_T = \int_t^{\t_k} \D Z_s dX_s^{t,x} + \sum_{i=0}^k [c(\d Z^2_{\t_i}) - c(\d Z^1_{\t_i})].
   $$
   Fix $i\le k$. If $\d Z^1_{\tau_i}\d Z^2_{\tau_i} >0$, then by Assumption ({\bf H3}) (iii) we get
   $$
   c(\d Z^2_{\t_i}) - c(\d Z^1_{\t_i})\le \rho(|\d Z^2_{\t_i}-\d Z^1_{\t_i}|) =\rho(|\D Z_{\t_{i-1}}-\D Z_{\t_i}|) \le \rho(z).
   $$
    Now assume $\d Z^1_{\tau_i}\d Z^2_{\tau_i}=0$. If $\d Z^1_{\tau_i}=0$, by Definition \ref{admissible} and (\ref{ti}) we must have  $i=0$ and $Z^1_{\t_0}=z$. This implies that $Z^2_{\t_0}=0$ and thus $\d Z^2_{\tau_0}=0$. If  $\d Z^1_{\tau_i}\neq 0$, then again we have  $\d Z^2_{\tau_i} =0$, and thus
   $$
   c(\d Z^2_{\t_i}) - c(\d Z^1_{\t_i})=- c(\d Z^1_{\t_i})\le 0\le \rho(z).
   $$
   Therefore, for some appropriately defined $\mathcal F_{T}$-measurable random 
   variable $\xi$, we have
  \beaa
  && E\{U(Y^{t,x,y,z,Z^1}_T)\}- V^n(t,x,y,0)\le E\Big\{U(Y^{t,x,y,z,Z^1}_T)- U(Y^{t,x,y,0,Z^2}_T)\Big\}\\
   && = E\Big\{U'(\xi)[Y^{t,x,y,z,Z^1}_T - Y^{t,x,y,0,Z^2}_T]\Big\}\le E\Big\{U'(\xi)\Big[\int_t^{\t_k}\D Z_s dX_s^{t,x} + k\rho(z)\Big]\Big\}\\
   && \le \L E\Big\{|\int_t^{\t_k}\D Z_s dX_s^{t,x}| + n\rho(z)\Big\}\le C[z +n\rho(z)].
 \eeaa
 This implies that
 $$
 V^n(t,x,y,z) - V^n(t,x,y,0)\le C[z +n\rho(z)].
 $$
 Sending $z\downarrow 0$ we obtain $V^n(t,x,y,0+)\le  V^n(t,x,y,0)$.

 Similarly, we can prove $V^n(t,x,y,0-)\le
 V^n(t,x,y,0)$. The proof is now complete.
 \qed

 \section{The Approximating Optimal Strategies}
  \setcounter{equation}{0}

In this section we construct the optimal strategy $Z^n\in
\cZ_t^n(z)$ for the approximating problem (\ref{Vn}).
We will provide the uniform estimate on
$Z^n$'s in next section.

We start with some auxiliary results. For any function $\f(t,x,y,z)$, define
\bea
 \label{fbardef}
 \left.\ba{lll}
 \bar\f(t,x,y,z) &:=& \sup_{\tilde z \in [-M,M]} \f\big(t,x,y-c(\tilde z - z), \tilde z\big);\\
  \hat\f(t,x,y,z) &:=& \sup_{\t\ge t} E\big[\bar \f(\t, X_\t^{t,x}, y+z[X_\t^{t,x}-x], z)\big],
  \ea\right.
 \eea
 where the supremum is taken over all stopping times $\t\ge t$. It is clear that
 \beaa
 \bar \f \le \hat \f &\mbox{and}& \hat\f(T,x,y,z) = \bar\f(T,x,y,z)
 \eeaa
 The following lemma is important for our construction of $Z^n$.
  \begin{lem}
 \label{fbar}
 Assume (H1)-(H3). Suppose that a function $\f: [0,T]\times\dbR^3\mapsto \dbR$ enjoys the
 following properties:

\ss
a) $|\f(t,x,y,z)|\le C[1+|y|]$;

\ss
b) $\f$ is increasing in $y$;  uniformly continuous in $(t,x,y)$; and uniformly continuous in $z$ in $[-M,0)$ and in $(0,M]$;

\ms
c) $\f(t,x,y,0+)\le \f(t,x,y,0)$, $\f(t,x,y,0-)\le \f(t,x,y,0)$.

\no Then

(i) $|\bar\f(t,x,y,z)|\le C[1+|y|]$ and $\bar\f$ is also uniformly continuous in $(t,x,y)$. Moreover,    there exists a Borel
measurable function $\psi(t,x,y,z)$ such that $|\psi|\le M$ and
 \bea
 \label{psi}
 \bar\f(t,x,y,z) = \f(t,x,y-c(\psi(t,x,y,z) - z), \psi(t,x,y,z)).
 \eea

 (ii) The optimal stopping problem $\hat \f$ admits an optimal stopping time $\t^*$:
 \beaa
 \t^* := \inf\Big\{s\ge t: \hat \f(s, X_s^{t,x}, y+z[X_s^{t,x}-x], z)=\bar \f(s, X_s^{t,x}, y+z[X_s^{t,x}-x], z)\Big\}.
 \eeaa
 \end{lem}

{\it Proof.} First, assume (i) holds true, then (ii) is a standard result in optimal stopping theory, see e.g.   \cite[Appendix
D]{KS}. To prove (i), note that
\beaa
|\bar\f(t,x,y,z)|\le C\sup_{\tilde z\in [-M,M]}\big[1+|y-c(\tilde z-z)|\big] \le C\big[1+ \sup_{\tilde z\in [-2M,2M]} |c(\tilde z)| + |y|\big]\le C[1+|y|].
\eeaa
Moreover,  by (H3) and the regularity of $\f$   we see that
$\f(t,x,y-c(\tilde z - z), \tilde z)$ is uniformly continuous in
$(t,x,y)$, uniformly in $(z,\tilde z)$. Thus  $\bar\f$ is uniformly
continuous in $(t,x,y)$.

It remains to construct the function $\psi$. We shall apply the measurable selection theorem in Wagner \cite{Wag77}.
For notational convenience, we define $\theta:= (t,x,y,z)\in [0,\infty)^2 \times \mathbb R \times [-M,M]$,
 $g(\theta,\tilde z) := \varphi(t,x,y - c(\tilde z - z), \tilde z)$, and $\bar g(\theta, \Gamma)
 := \sup_{\tilde z \in \Gamma} g(\theta, \tilde z)$ for any  Borel set $\Gamma \subset [-M, M]$ (by convention $\bar g(\theta, \emptyset) := -\infty$).
Consider a set-valued function defined by
$$F(\theta) = \{ z'\in [-M, M]: g(\theta, z') = \sup_{\tilde z\in [-M, M]}g(\theta, \tilde z)\}.
$$
By our conditions, one may easily check that $g$ is upper semicontinuous in $\tilde z$. Then  $F(\theta)$ is a nonempty and  closed set for  any $\theta$ in the domain
 $ [0,\infty)^2\times \mathbb R \times [-M,M]$. In light of \cite[Theorem 4.1]{Wag77}, to obtain the measurable $\psi$ it suffices to prove: 
\bea
\label{MST}
\mbox{for any open set $\Gamma\subset [-M, M]$,
 $\{\theta: F(\theta) \cap \Gamma \neq \emptyset\} \subset R^{4}$ is a Borel set.  
 }
 \eea
To see this, we first assume $c(\cdot)$ is continuous.
 Since $\varphi(t,x,y,z)$ is continuous in $ [0,\infty)^2 \times \mathbb R \times [-M,0)$, $g(\theta, \tilde z)$ is also continuous in 
 $ [0,\infty)^2 \times \mathbb R \times [-M,0) \times [-M,0)$.
 Therefore, if $\Gamma \subset [-M,0)$, then we can write,
 denoting the set of all rational numbers by $\mathbb Q$, that
 $$
\bar\varphi(\th):= \bar g(\theta, \Gamma)
 = \sup_{\tilde z \in \Gamma} g(\theta, \tilde z)
 = \sup_{\tilde z \in \Gamma \cap \mathbb Q} g(\theta, \tilde z)
 $$
 Thus, $\bar g(\cdot,\Gamma)$ is a Borel measurable function (in fact, it is a Baire function of Class 1) for 
  $\Gamma \subset [-M,0)$. Similar argument shows that 
 $\bar g(\cdot, \Gamma)$ is also Borel measurable if
 $\Gamma\subset (0, M]$. On the other hand,
 if $\Gamma = \{0\}$, then $\bar g(\theta, \Gamma) = g(\theta,0) 
 = \varphi (t,x,y-c(-z),0)$ is obviously continuous.
 In general, if $\Gamma \subset [-M, M]$ is an open set, 
 we can partition this set into $\Gamma = \cup_{i = 1,2,3} \Gamma_{i}$,
 where 
 $\Gamma_{1} = \Gamma \cap [-M,0)$, 
 $\Gamma_{2} = \Gamma \cap (0, M]$,  and
 $\Gamma_{3} = \Gamma \cap \{0\}$. Then, we can see
 $\bar g(\cdot, \Gamma)$ is Borel measurable, since
  $\bar g(\theta, \Gamma) = \max_{i = 1,2,3} \bar g(\theta, \Gamma_{i})$.
 Therefore, noting that $\bar g(\th,\G)\le \bar g(\th, [-M,M])$ as $\G\subset [M,M]$, we can conclude that the set 
 $$\{\theta: F(\theta) \cap \Gamma \neq \emptyset\} = \{\theta: \bar g(\theta, \Gamma) = \bar g(\theta, [-M, M])\},$$
 whence a Borel set, and thus (\ref{MST}) holds when $c$ is continuous at $0$. In the general case,  since $c$ is lower semicontinuous
 at $0$, one can prove (\ref{MST}) by repeating the above arguments but with 
 the utilization of 
 $\bar g(\theta, \Gamma) = 
 \max\{\sup_{\tilde z \in \Gamma \cap \mathbb Q} g(\theta, \tilde z), g(\theta, z)\}
 $. 
 \qed

We now give the main existence result of $Z^n$ for this section.
\begin{thm}
\label{Znexist}
Assume (H1)--(H3). Then, for each $n$ and any fixed $(t,x,y,z)$,
 \bea
 \label{Vnn-1}
 V^n(t,x,y,z) = \hat V^{n-1}(t,x,y,z)
 \eea
Moreover, there exists an optimal $Z^n\in\cZ^n_t(z)$ such that $V^n(t,x,y,z) = E\Big[U(Y^{t,x,y,z,Z^n}_T)\Big]$.
\end{thm}
 {\it Proof.}  We proceed in several steps.

 {\it Step 1.} We first show that
  \bea
 \label{Vnn-1hat}
 V^n(t,x,y,z) \le \hat V^{n-1}(t,x,y,z)
 \eea
 Indeed, let $Z\in \cZ^n_t(z)$. If $Z_{\t_0} \neq z$, then $Z\in \cZ^{n-1}_t(Z_{\t_0})$, and
\beaa
 E\Big[U(Y^{t,x,y,z,Z}_T)\Big] &=& E\Big[U(Y^{t,x,y-c(Z_{\t_0}-z), Z_{\t_0}, Z}_T)\Big] \\
&\le& V^{n-1}(t,x,y-c(Z_{\t_0}-z), Z_{\t_0}) \le \bar V^{n-1}(t,x,y,z)\le  \hat V^{n-1}(t,x,y,z).
\eeaa
If $Z_{\t_0}=z$, then we denote
\bea
\label{Y0}
X_s:= X_s^{t,x} &\mbox{and}& Y^0_s:= y+ z[X_s-x],\qq s\in [t,T].
\eea
Clearly we have  $Z\in \cZ^{n-1}_{\t_1}(Z_{\t_1})$ and
 \beaa
  Y^{t,x,y,z,Z}_T&=& y+z(X^{t,x}_{\t_1}-x)-c(Z_{\t_1}-z)+\int_{\t_1}^T Z_s
  dX^{t,x}_s-
  \sum_{i=2}^\infty c(\d Z_{\t_i})\\
  &=&  Y^{\t_1, X_{\t_1}, Y^0_{\t_1}-c(Z_{\t_1}-z), Z_{\t_1},Z}_T.
\eeaa
This implies that
 \beaa
 \label{barV}
 E\Big[U(Y^{t,x,y,Z}_T)\Big] &\le& E\Big[ V^{n-1}(\t_1,X_{\t_1},Y^0_{\t_1}-c(Z_{\t_1}-z), Z_{\t_1})\Big]\\
 &\le& E\Big[\bar V^{n-1}(\t_1,X_{\t_1},Y^0_{\t_1}, z)\Big]\le \hat V^{n-1}(t,x,y,z).
 \eeaa
Since $Z$ is arbitrary, we obtain (\ref{Vnn-1hat}).

{\it Step 2.} We now construct $Z^n$. By the results in Section 3, we see that  we may apply Lemma \ref{fbar}  on $\f := V^{n-1}$.  Let $\psi$ and $\t^n_1:=\t^*$ be given as in Lemma \ref{fbar} (ii).  Set
$Z^n_s := z$, for $s\in [t, \t^n_1)$, and $Z^n_{\t^n_1} := \psi(\t^n_1, X_{\t^n_1}, Y^0_{\t^n_1}, z)$. Then by Lemma  \ref{fbar}  
we get
\bea
\label{hatV0}
\hat V^{n-1}(t,x,y,z) = E\Big[V^{n-1}(\t^n_1, X_{\t^n_1}, Y^0_{\t^n_1} - c(Z^n_{\t^n_1}-z),Z^n_{\t^n_1})\Big].
\eea
We remark that if $\t^n_1 =t$, then $Z^n$ has a jump at $t$, and if $\t^n_1>t$, then $Z^n_t = z$ and does not jump at $t$.  Note that $Y^0_{\t^n_1} - c(Z^n_{\t^n_1}-z) = Y^{t,x,y,z, Z^n}_{\t^n_1}$. Then, by  (\ref{Vnn-1hat}) we obtain
\bea
\label{hatV1}
V^n(t,x,y,z) \le E\Big[V^{n-1}(\t^n_1, X_{\t^n_1}, Y^{t,x,y,z, Z^n}_{\t^n_1},Z^n_{\t^n_1})\Big].
\eea
Repeating the above arguments, we define $\t^n_i$, $i=2,\cds, n-1$ and $Z^n$ on $[t, \t^n_{n-1}]$ such that
\bea
\label{hatV2}
V^{n-i+1}(\t^n_{i-1}, X_{\t^n_{i-1}}, Y^{t,x,y,z, Z^n}_{\t^n_{i-1}}, Z^n_{\t^n_{i-1}}) \le E_{\t^n_{i-1}}\Big[V^{n-i}(\t^n_i, X_{\t^n_i}, Y^{t,x,y,z,Z^n}_{\t^n_i},Z^n_{\t^n_i})\Big].
\eea
Finally, for $V^1$, there exists $\t^n_n \ge \t^n_{n-1}$ such that, by setting $Z^n_s:= Z^n_{\t^n_{n-1}}$ for $s\in [\t^n_{n-1}, \t^n_n)$ and  $Z^n_s:= Z^n_{\t^n_{n-1}}$ for $s\in[\t^n_n,T]$,
\bea
\label{hatV3}
V^{1}(\t^n_{n-1}, X_{\t^n_{n-1}}, Y^{t,x,y,z,Z^n}_{\t^n_{n-1}},Z^n_{\t^n_{n-1}}) = E_{\t^n_{n-1}}\Big[U(Y^{t,x,y,z,Z^n}_{\t^n_n})\Big] =  E_{\t^n_{n-1}}\Big[U(Y^{t,x,y,z,Z^n}_T)\Big].
\eea

Now combining (\ref{hatV1})-(\ref{hatV3}) we obtain
\beaa
V^n(t,x,y,z) \le  E\Big[U(Y^{t,x,y,z,Z^n}_T)\Big].
\eeaa
Since clearly $Z^n \in\cZ^n_t(z)$,  it is an optimal strategy for the optimization problem $V^n$.

{\it Step 3.} Since $V^n(t,x,y,z) =  E\Big[U(Y^{t,x,y,z,Z^n}_T)\Big]$. By Step 2 we see that (\ref{hatV1}) (and (\ref{hatV2})) should hold with equality. This, together with (\ref{hatV0}), implies (\ref{Vnn-1}).
 \qed

\section{Regularity of the Value Function}
\setcounter{equation}{0}

In this section we give some {\it uniform estimates} of the value
function $V$. We should note that the regularity of $V$ with
respect to the variables $(t,x,y)$ are clear, since the estimates
(\ref{regx}), (\ref{regy}), and (\ref{regt}) in Proposition
\ref{regxy} are already uniform with respect to $n$. The estimate
(\ref{Vnregz}), however, depends heavily on $n$. In fact, in the case
$|z|=|z|^\a$, $0<\a<1$, one can check that
$\rho_n(|z|)=n^{1-\a}|z|^\a\to\infty$. Therefore the regularity of
$V$ with respect to $z$ is by no means clear.

We first take a closer look at the approximating optimal strategies
$\{Z^n\}_{n=1}^\infty$. Since our purpose is to construct the
optimal piecewise constant control,
it is thus conceivable that a uniform bound on $N(Z^n)$ would 
be extremely
helpful.

We begin by considering the case where a fixed cost is present.
For each $(t,x,y,z)$, we denote $Z^n$ to be the optimal portfolio
for $V^n(t,x,y,z)$, when the context is clear.

\begin{prop}
\label{EN1} Assume (H1)--(H3), and assume further that $c(z)\ge
c_0>0$ for any $z\neq 0$.
Then there exists a constant $C>0$ such that
 \bea
 \label{boundN}
 E\{N(Z^n)\}\le {C\over \l c_0}, \q \mbox{\rm for all $n$ and all
 $(t,x,y,z)$.}
 \eea
\end{prop}
{\it Proof.} Denote $Z^0:= z\one_{[t,T)}\in \cZ^1(z)$. Then
$$
Y^{t,x,y,z,Z^0}_T- Y^{t,x,y,Z^n}_T = \sum_{i=0}^n c(\d Z^n_{\t_{i}})+\int_t^T [z-Z^n_s]dX^{t,x}_s  - c(-z).
$$
Note that $V^n$'s are non-decreasing in $n$. Then
 \beaa
0&\ge& V^0(t,x,y,z) - V^n(t,x,y,z)\ge E\{U(Y^{t,x,y,z,Z^0}_T)\}-E\{U(Y^{t,x,y,Z^n}_T)\}\\
&\ge& \l E\Big\{\sum_{i=0}^n c(\d Z^n_{\t_{i}})\Big\}
-\L E\Big\{\Big|\int_t^T [z-Z^n_s]dX^{t,x}_s\Big| + c(-z)\Big\}\ge \l c_0 E\{N(Z^n)\} - C.
 \eeaa
The result follows immediately. \qed

We next investigate the problem under (H4).  We first have the following technical lemma.
\begin{lem}
\label{lem-concave}
Assume (H1)-(H4) hold. Denote:
\bea
\label{C01}
\left.\ba{c}
\dis\a_1 := { 1-\eta_2\over \eta_2},\q  \b_1 := {1-\eta_2\over 1+\g};\\
\dis C_0 := {\L\over\l}\Big[\|b\|_\infty T + \|\si\|_\infty\sqrt{T} +L_0\Big]+1;\q
C_1:= C_0\Big[2+\L({1\over\a_1}+{1\over \b_1})\Big],
\ea\right.
\eea
There exists  a constant $\e_1 \in (0,{\e_0}]$ such that, for any $0<|z_1|< \e_1$,

(i)  $c(z_1)\ge C_0|z_1|$.

(ii) For any $z_2 \ge z_1>0$ or $z_2\le z_1<0$,  we have
\beaa
 c(z_1)+c(z_2) - c(z_1+z_2) \ge \Big[\a_1[c(z_1+z_2)-c(z_2)]\Big] \vee \Big[ \b_1[c(-z_1-z_2)-c(-z_2)]\Big].
\eeaa

\ss (iii) For any $z_2 \ge {1\over 2}z_1 >0$, or $z_2\le {1\over 2}z_1 <0$, or $|z_2|>|z_1|$,   we have
\beaa
c(z_1)+c(z_2)-c(z_1+z_2)\ge C_1|z_1|.
\eeaa
\end{lem}
{\it Proof.}  For $\th = {3\over 2}, 2, 3$, set
\beaa
\tilde\eta_\th :=  {1\over 2}[\eta_\th+ 1], &\mbox{so that}& \eta_\th<\tilde\eta_\th < 1.
\eeaa

(i)  By (\ref{c0infty}), there exists $0<\e\le {\e_0}$ such that $c(2z) \le 2\tilde\eta_2 c(z)$ for all $|z|\le \e$. By induction one can easily show that ${c(2^{-n}\e)\over 2^{-n}\e} \ge {c(\e)\over \e\tilde\eta_2^{n}}$. Fix $n_0$ such that ${c(\e)\over \e\tilde\eta_2^{n_0}} \ge 2C_0$, and set $\e_1 := 2^{1-n_0}\e$. For for any $0<|z|< \e_1$, there exists $n\ge n_0$ such that $2^{-n}\e < |z|\le 2^{1-n}\e$. Then
\beaa
{c(z) \over |z|} \ge {c(2^{-n}\e)\over 2^{1-n}\e} \ge {1\over 2} {c(\e)\over \e\tilde\eta_2^{n}} \ge {1\over 2} {c(\e)\over \e\tilde\eta_2^{n_2}} \ge C_0.
\eeaa

(ii) Without loss of generality, we assume $z_2\ge z_1>0$.  We may rewrite the required inequality as
\beaa
&f(z_1,z_2) \le c(z_1) ~\mbox{where}~&\\
&f(z_1,z_2):=[c(z_1+z_2)-c(z_2)] +  \Big[ \a_1[c(z_1+z_2)-c(z_2)]\Big]\vee \Big[  \b_1[c(-z_1-z_2)-c(-z_2)]\Big].&
\eeaa
If $z_2 \in [z_1, \e_0]$, by the concavity of $c$,  $f(z_1, z_2)$ is decreasing in $z_2$,  then
\beaa
f(z_1, z_2) \le f(z_1, z_1) = [c(2z_1)-c(z_1)] +  \Big[ \a_1[c(2z_1)-c(z_1)]\Big]\vee \Big[  \b_1[c(-2z_1)-c(-z_1)]\Big].
\eeaa
By choosing $\e_1$ small enough, we have
\beaa
c(2z_1) -c(z_1)\le [2\tilde\eta_2 -1]c(z_1) = \eta_2c(z_1) &\mbox{and}& c(-2z_1)-c(-z_1)\le (1+\g) c(z_1).
\eeaa
Then
\beaa
f(z_1, z_1) \le \Big[\eta_2 + [(\a_1 \eta_2)\vee (\b_1 (1+\g))]\Big]c(z_1) =c(z_1).
\eeaa
If  $z_2 \in [\e_0, 2M]$, by (H4)-(ii) we have
\beaa
f(z_1, z_2) \le L_0 z_1 +  [ \a_1\vee\b_1] L_0 z_1 =[ 1+  \a_1\vee\b_1] L_0 z_1.
\eeaa
By replacing $C_0$ with $[ 1+  \a_1\vee\b_1] L_0$ and setting $\e_1$ smaller if necessary, it follows from  (i) that $f(z_1, z_2) \le c(z_1)$.

(iii)  Without loss of generality, we assume $z_1>0$, and it suffices to show that
\beaa
g(z_1,z_2):=c(z_1+z_2)-c(z_2) + C_1|z_1| \le c(z_1).
\eeaa
 If $z_2 \le - z_1$, then $z_2 < z_1+z_2 \le 0$, and thus $g(z_1,z_2) \le C_1|z_1|$.
By setting $\e_1$ smaller if necessary,  the result follows from the proof of (i) by replacing $C_0$ with $C_1$.

If $z_2 \ge \e_0$, then $g(z_1, z_2) \le  [L_0+C_1]|z_1|$. The result follows from the proof of (i) by replacing $C_0$ with $L_0+C_1$.
Finally, if ${1\over 2}z_1\le z_2 \le \e_0$, then $g(z_1,z_2)$ is decreasing in $z_2$, and thus
 \beaa
&&g(z_1,z_2) \le g(z_1, {1\over 2}z_1) =  c({3z_1\over 2}) - c({z_1\over 2}) + C_1z_1.
\eeaa
Then, by choosing $\e_1$ smaller if necessary, we have
\beaa
c(z_1) - g(z_1, z_2) &\ge&    [c(z_1) - {2\over 3}c({3z_1\over 2})] +[ c({z_1\over 2})-{1\over 3}c({3z_1\over 2})] - C_1z_1\\
 &\ge& [1-\tilde \eta_{3\over 2}] c(z_1) + [1-\tilde \eta_{3}] c({z_1\over 2}) -C_1z_1,
 \eeaa
 Now the result follows from the proof of (i) by replacing $C_0$ with an appropriate larger constant.
\qed

 \ms

To extend Proposition \ref{EN1} under (H4), we need an analysis on the number of the small jumps. For
this purpose, we fix the constants  $\e_1$, $C_0$, and $C_1$  given in Lemma \ref{lem-concave}.  Define:
 \be
 \label{Ai}
 A^n_i := \{0<|\d Z^n_{\t_i}|<\e_1\},\q B^n_i :=
 \{|\d Z^n_{\t_i}|\ge \e_1\}, \q i=0,\cds, n; ~ n>0,
 \ee
 The following result is crucial.
\begin{thm}
\label{smalljump} Assume (H1)--(H4). Then for any fixed $m$,
 \bea
 \label{pAi}
 P\Big(\sum_{i=0}^n\one_{A^n_i} \ge m\Big) \le {1\over
 2^{m}},\q\forall n\ge m.
 \eea
\end{thm}

The proof of Theorem \ref{smalljump} depends heavily on the
following technical result, whose proof is quite lengthy and will be
deferred to Section 7 in order not to distract the discussion.
\begin{prop}
\label{smalllarge} Assume (H1)--(H4).
Then, for any $n$ and $i$,
 $P$-a.s in $A^n_i$ one has: (i)
$P\{B^n_{i+1}|\cF_{\t_i}\} \le {C_0\over C_1}<\frac12$
for the constants $C_{0}$ and $C_{1}$ defined in \eqref{C01}, 
and (ii) $Z^n_{\t_i}=0$.
\end{prop}
[{\it Proof of Theorem \ref{smalljump}}.] Define $k_{-1}:= -1$, and
$$
k_{j}:=  \inf\{i>k_{j-1}: 0<|\d Z^n_{\t_i}|<\e_1\}\wedge (n+1),\q j=0,1,\cds,n.
 $$
 Then
 $ P\Big(\sum_{i=0}^n\one_{A^n_i} \ge m\Big) = P(k_m\le n)$.
We claim that, for each $0\le j<n$,
 \bea
 \label{kjlen}
 \{k_{j+1}\le n\} \subseteq A^n_{k_j}\cap B^n_{k_j+1}, \q \mbox{$P$-a.s.}
 \eea
(It is important to note here that the left side contains $k_{j+1}$ while the superscript of $B$ on the right side is $k_j+1$
!)

Indeed, we first note that $\{k_{j+1}\le
n\}\subset\{k_j\le n\} \subseteq A^n_{k_j}$, and consider the set
$A^n_{k_j}\setminus B^n_{k_j+1}$.
Suppose that $Z^n_{\t_{k_j+1}}\neq Z^n_{\t_{k_j}}$ on
$A^n_{k_j}\setminus B^n_{k_j+1}$. Then
$0<|Z^n_{\t_{k_j+1}}-Z^n_{\t_{k_j}}|<\e_1$, and by Proposition
\ref{smalllarge} (ii) we must have both $Z^n_{\t_{k_j}}=0$ and
$Z^n_{\t_{k_j+1}}=0$, $P$-a.s., a contradiction.
Thus we must have $Z^n_{\t_{k_j+1}}=Z^n_{\t_{k_j}}$ on
$A^n_{k_j}\setminus B^n_{k_j+1}$. Then by the definition of $\t_i$ in (\ref{ti})
we know $\t_{k_j+1}=T$ and thus
$Z^n_{\t_{k_j}}=Z^n_{\t_{k_j+1}}=\cds=Z^n_{\t_n} = 0$. Namely
$k_{j+1}=n+1$. In other words, $A^n_{k_j}\setminus B^n_{k_j+1}\subseteq
\{k_{j+1}=n+1\}$. Note that $\{k_{j+1}\le n\}\subseteq
A^n_{k_j}\setminus \{k_{j+1}=n+1\}$, (\ref{kjlen}) follows.

  Next, applying Proposition \ref{smalllarge} (i) we derive from (\ref{kjlen}) that
\beaa
&&P\big(\sum_i \one_{A^n_i} \ge m\big) = P(k_m \le n)\le P\Big(\bigcap_{j=0}^{m-1}[A^n_{k_j}\cap B^n_{k_j+1}]\Big)\\
&& = E\Big\{\prod_{j=0}^{m-1} [\one_{A^n_{k_j}}\one_{B^n_{k_j+1}}]\Big\}= E\Big\{\Big[\prod_{j=0}^{m-2}[ \one_{A^n_{k_j}}\one_{B^n_{k_j+1}}]\Big]\one_{A^n_{k_{m-1}}}E\{\one_{B^n_{k_{m-1}+1}}|\cF_{\t_{k_{m-1}}}\}\Big\}\\
&&\le E\Big\{\Big[\prod_{j=0}^{m-2}[
\one_{A^n_{k_j}}\one_{B^n_{k_j+1}}]\Big]\one_{A^n_{k_{m-1}}}{1\over 2}\Big\}\le
{1\over 2}E\Big\{\prod_{j=0}^{m-2}[ \one_{A^n_{k_j}}\one_{B^n_{k_j+1}}]\Big\}.
\eeaa
 Repeating the argument $m-1$ more times we prove the theorem. \qed

The following theorem is a generalized version of Proposition
\ref{EN1}.
\begin{thm}
\label{EN2} Assume Assumptions (H1)--(H4). Then it holds that
$$
E\{N(Z^n)\} \le C\Big[1+{1\over c(\e_1)\wedge c(-\e_1)}\Big]<\infty,
\q\forall n.
$$
\end{thm}
{\it Proof.}  Denote
$$N_1(Z^n):=
\sum_{i=0}^n \one_{A^n_i}, \q
N_2(Z^n) := \sum_{i=0}^n
\one_{B^n_i}.
 $$
Then $E\{N(Z^n)\} = E\{N_1(Z^n)\} + E\{N_2(Z^n)\}$. First, Theorem \ref{smalljump} implies that
 \be
 \label{ENest}
 E\{N_1(Z^n)\} = \sum_{m=0}^n P(N_1(Z^n)\ge m) \le \sum_{m=0}^n
 {1\over 2^{m}} \le 2,
 \ee
 Next, one can estimate $E\{N_2(Z^n)\}$ along the lines as Proposition
\ref{EN1}. Indeed, note that
 \beaa
&& E\Big\{U(y+\int_t^T Z^n_s dX_s)\Big\} - V^n(t,x,y,z)\\
 &=& E\Big\{U(y+\int_t^T Z^n_s dX_s)-U(y+\int_t^T Z^n_s dX_s-\sum_i c(\d Z^n_{\t_i}))\Big\}\\
& \ge& \l E\Big\{\sum_i c(\d Z^n_{\t_i})\Big\} \ge \l
E\Big\{\sum_i (c(\e_1)\wedge
c(-\e_1))\one_{B^n_i}\Big\}= \l [c(\e_1)\wedge c(-\e_1)] E\{N_2(Z^n)\}.
 \eeaa
On the other hand, recalling (\ref{V1}) we have
 \beaa
 && E\{U(y+\int_t^T Z^n_s dX_s)\} - V^n(t,x,y,z) \le E\{U(y+\int_t^T Z^n_s dX_s)\}
 - V^1(t,x,y,z)\\
 &=&\sup_{\t\ge t}\Big|E\Big\{U(y+\int_t^T Z^n_s dX_s)-U(y+\int_t^\t z dX_s -
 c(-z))\Big\}\Big|\\
 &\le& \L E\Big\{|\int_t^\t (Z^n_s-z)dX_s|+| \int_\t^T Z^n_sdX_s|+ c(-z)\Big\}\le
 C\L,
 \eeaa
Then
$
E\{N_2(Z^n)\}\le {C\L\over\l [c(\e_1)\wedge c(-\e_1)]}.
$
This, together with (\ref{ENest}), proves the theorem.
%
\qed

As  a consequence Theorem \ref{EN2}, we have the
second main result of this section, which improves (\ref{Vnregz}) and whose proof is also postponed to Section 7.
\begin{thm}
\label{regzunif} Assume (H1)--(H3). Assume further that either
$c(z)\ge c_0>0$, for all $z\neq 0$ or (H4) holds. Then there exists
a generic constant $C>0$, such that for any $z_1, z_2$ with the same
sign, and for all $n$, it holds that
 \bea
 \label{uniregz}
  &&|V^n(t,x,y,z_1)-V^n(t,x,y,z_2)|  \le C[|\D z|+\rho(|\D z|)];\\
 \label{uniconv}
&& V(t,x,y,z)-V^n(t,x,y,z)\le {C\over n}.
 \eea
\end{thm}

As the direct consequences of Propositions \ref{regxy} and \ref{regz0},
 and Theorem \ref{regzunif} we have
\begin{thm}
\label{Vreg} Assume(H1)--(H3), and assume either $c(z)\ge c_0$,
$z\neq 0$ or (H4). Then

(i) $|V(t,x_1,y,z)-V(t,x_2,y,z)|\le C|\D x|$.

(ii) $\l \D y\le V(t,x,y_1,z)-V(t,x,y_2,z)\le \L \D y, \forall \D y:= y_1-y_2\ge 0$.

(i) $|V(t_1,x,y,z)-V(t_2,x,y,z)|\le C|\D t|^{1\over 2}$.

(iv) $|V(t,x,y,z_1)-V(t,x,y,z_2)|\le C[|\D z|+\rho(|\D z|)], \forall z_1, z_2$ with the same sign.

(v) $V(t,x,y, 0+)\le V(t,x,y,0), V(t,x,y,0-)\le V(t,x,y,0)$.
\end{thm}

\section{The Optimal Strategy $Z^*$}
\setcounter{equation}{0}

In this section we construct the optimal controls for the
original problem (\ref{V}).
We should note that by virtue of Proposition \ref{EN1} and Theorem \ref{EN2},
one can easily show that under our assumptions $Z^n$ should converge
in distribution. But this does not seem to be helpful for our
construction of the optimal strategy. In fact, in general we will
have to extend the probability space, and it is not clear whether
the limit process will have the desired adaptedness that is
essential in our application. We thus construct the optimal portfolio $Z^*$ for (\ref{V}) directly.
%

In light of
the construction of the optimal strategy $Z^n$, we know that the
function $\bar V=V$ should play the role of an ``obstacle" that will
trigger the jumps, as it is usually the case in impulse control
literature.
%
To this end let us consider the following set
 \bea
 \label{cO}
 \cO(z) := \{(t,x,y): V(t,x,y,z) > V(t,x, y-c(\tilde z-z),\tilde z), \forall
 \tilde z\neq z\},~~\cO:= \bigcup_z \cO(z).
 \eea
Intuitively, the set $\cO(z)$ should define an``inaction region",
since a change of position (on $z$) would decrease the value
function. Furthermore, following the standard impulse control theory
one would expect that $\cO(z)$ is an open set and the trade will
take place when $(t,x,y)\in\pa\cO(z)$. This is indeed the case when $c(z)\ge c_0>0$ for $z\neq 0$. However, unfortunately in our more general case we
only have the following result.
\begin{lem}
\label{On}
 Assume (H1)--(H4). Define
 \bea
 \label{cOn}
 \cO_n(z):= \{(t,x,y): V(t,x,y,z) > V(t,x, y-c(\tilde z-z),\tilde z), \forall |\tilde z - z|\ge {1\over n}\}.
 \eea
Then  $\cO_n(z)$ is open, for all $n$, and $\cO(z)=\bigcap_n
\cO_n(z)$.
\end{lem}

{\it Proof.} Denote \be \label{Vn0} \bar V_n(t,x,y,z) := \sup_{|\tilde
z-z|\ge {1\over n}}V(t,x, y-c(\tilde z-z),\tilde z). \ee Apply
Theorem \ref{Vreg} and follow the proof of Lemma \ref{fbar}, we
know $\bar V_n$ is continuous in $(t,x,y)$ and there exists a Borel
measurable function $\psi_n$ such that $|\psi_n(t,x,y,z)-z|\ge
{1\over n}$ and
 $$
 V(t,x, y-c(\psi_n(t,x,y,z)-z),\psi_n(t,x,y,z))=\bar V_n(t,x,y,z).
 $$
 This implies that
  $$
\cO_n(z) = \{(t,x,y): V(t,x,y,z)>\bar V_n(t,x,y,z)\}
$$
and thus $\cO_n(z)$ is open. That $\cO(z)=\cap_{n=1}^\infty\cO_n(z)$
is obvious. The proof is complete.
\qed

We remark that Lemma \ref{On} does not imply that the set $\cO(z)$
is an open set. Therefore, if we follow the scheme in the previous
sections to define, for given $(t,x,y,z)\in\cO$ and recalling (\ref{Y0}),
 \bea
 \label{tau}
 \t := \inf\{s\ge t: (s,X_s, Y^0_s)\notin \cO(z)\}\wedge T.
 \eea
Then intuitively it is possible that $P\{\t=t\}>0$ and/or
$P\{(\t,X_\t,Y^0_\t)\in\cO(z)\}>0$. In either case the construction
procedure will fail. The following Theorem, which excludes the above cases, is therefore essential.
 \begin{thm}
 \label{tnlimit}
 Assume (H1)--(H4). Define, for each $(t,x,y,z)\in\cO$ and $n>0$,
 \bea
 \label{taun}
 \t^n:= \inf\{s\ge t: (s,X_s, Y^0_s)\notin \cO_n(z)\}\wedge T,
 \eea
 and let $\t$ be defined by (\ref{tau}).
 Then

 (i) $\t^n$ are decreasing stopping times and
  $(\t^n, X_{\t^n}, Y^0_{\t^n})\notin \cO_n(z)$ whenever $\t^n<T$.

 (ii) $\t^n\downarrow \t$ and thus $\t$ is also a stopping time.

 (iii) $P(\t^n>\t,\forall n)=0$
 and thus, $P$-a.s., $(\t, X_{\t}, Y^0_{\t})\notin \cO(z)$ when $\t<T$. In particular, this implies that $\t>t$.

 (iv) $V(t,x,y,z) = E\{V(\t, X_\t, Y^0_\t, z)\}.$
 \end{thm}

The proof of Theorem \ref{tnlimit} will depend heavily on an
important, albeit technical, lemma that characterizes the possible
behavior of the small jumps under our basic assumptions on the
liquidity/transaction cost function. The proof of this lemma is again rather
tedious, and we defer it to Section 7.
 \begin{lem}
 \label{cOcompact}
 Assume (H1)--(H4) and let $\e_1$ be that in Lemma \ref{lem-concave}. Suppose that for given $(t,x,y,z)$, $\tilde z$ is such that
 $0<|\tilde z-z|<\e_1$ and
 $V(t,x,y,z) = V(t,x,y-c(\tilde z-z),\tilde z)$, 
 then $\tilde z=0$.
\end{lem}

 [{\it Proof of Theorem \ref{tnlimit}}] (i) That $\t^n$'s are
 decreasing stopping times is obvious by definition. Also, since
 each $\cO_n(z)$ is an open set, thanks to Lemma \ref{On}, it follows
 immediately that  $(\t^n, X_{\t^n}, Y^0_{\t^n})\notin \cO_n(z)$, whenever
 $\t^n<T$.

 (ii) Denote $\dis\t^\infty := \lim_{n\to\infty}\t^n$. Since $\cO_n \supseteq  \cO$,
we have $\t^n\ge \t$ for any $n$ and thus $\t^\infty\ge \t$,
$P$-a.s. The claim is trivial when $\t=T$. Now assume $\t(\o)<T$. Then for any $\e>0$, there exists
$s<\t(\o)+\e$ such that $(s,X_s, Y^0_s)\notin \cO(z)$. Since
$\cO(z)=\bigcap_n \cO_n(z)$, there exists $n:=n(\o)$ such that
$(s,X_s(\o), Y^0_s(\o))\notin \cO_n(z)$. Thus $\t^n(\o)\le
s<\t(\o)+\e$ and therefore $\t^\infty(\o)< \t(\o)+\e$. Since $\e$ is
arbitrary, we get $\t^\infty\le \t$, and hence $\t^\infty=\t$.

 (iii) Choose $n_0$ such that $n_0>\max\{ \frac{1}{\e_1}, \frac1{|z|}
{\bf 1}_{\{z\neq 0\}}\}$, and note
that $\{\t^n>\t,\forall n\}\subset \{\t<T\}$. On $\{\t<T\}$, for
$n\ge n_0$ large enough, by (ii) we have $\t^n<T$ and thus there
exists $Z_{\t^n}$ such that $|Z_{\t^{n}}-z|\ge {1\over n}$ and
 $
 V(\t^{n}, X_{\t^{n}}, Y^0_{\t^{n}},z) = V(\t^{n}, X_{\t^{n}}, Y^0_{\t^{n}}
 -c(Z_{\t^{n}}-z),Z_{\t^{n}}).
 $
By Lemma \ref{cOcompact}, either $Z_{\t^{n}}=0$ or
$|Z_{\t^{n}}-z|\ge \e_1$. If $z=0$, then $Z_{\t^{n}}\neq 0$ and thus
$|Z_{\t^{n}}-z|\ge \e_1\ge {1\over n_0}$. If $z\neq 0$, then either
$|Z_{\t^{n}}-z|=|z|\ge {1\over n_0}$ or  $|Z_{\t^{n}}-z|\ge \e_1\ge
{1\over n_0}$. So in all the cases we have  $|Z_{\t^{n}}-z|\ge
{1\over n_0}$. This implies that $\t^n=\t^{n_0}$ for all $n$ large
enough. Therefore, $\t=\t^{n_0}$ and thus $(\t, X_{\t},
Y^0_{\t})\notin \cO(z)$.

 (iv) We first note that, taking $\t$ as the first trading time,
 we should have
 $$
 E\{V(\t, X_\t, Y^0_\t, z)\} = \sup\{E\{U(Y^{t,x,y,z,Z}_T)\}: Z\in \cZ_t, Z_s = z ~~{\rm for}~~ \forall s< \t\}.
 $$
It then follows that $E\{V(\t, X_\t, Y^0_\t, z)\}\le V(t,x,y,z)$.

On the other hand, note that $\bF$ is
quasi-left continuous, we can choose a sequence of stopping
times $\t_m\uparrow \t$ such that $\t_m<\t$ whenever $\t>t$.  We
claim that
 \be
 \label{Vm}
 V(t,x,y,z)\le E\Big\{V(\t_m, X_{\t_m}, Y^0_{\t_m}, z) \Big\}.
 \ee
Then by sending $m\to\infty$ we prove the theorem.

To prove (\ref{Vm}), we recall (\ref{Vn0}) and choose $n_0$ as in
(iii). On the set $\{\t>t\}$ and for $t\le s<\t$, denote
$$
I_s := V(s, X_s, Y^0_s, z) - \bar V_{n_0}(s,X_s,Y^0_s,z).
$$
By the proof of Lemma \ref{On} we have $I_s >0$. Since $I$ is continuous in $s$, we get
\be
\label{Is}
I^m:= \inf_{s\le \t_m} I_s >0.
\ee

For any $n\ge n_0$, let $Z^n$ be the optimal portfolio of
$V^n(t,x,y,z)$. If $Z^n_t \neq z$, by
Proposition \ref{smalllarge} (ii) and following similar arguments as in (iii),
we have $|Z^n_t-z|\ge {1\over n_0}$. Then
$$
V^n(t,x,y,z) = V^{n-1}(t,x,y-c(Z^n_t-z), Z^n_t) \le V(t,x,y-c(Z^n_t-z),Z^n_t) \le \bar V_{n_0}(t,x,y,z).
 $$
Thus, by (\ref{uniconv}),
$$
V(t,x,y,z) \le V^n(t,x,y,z)+{C\over n} \le \bar V_{n_0}(t,x,y,z)+{C\over n},
$$
and therefore $Z^n_t=z$ for $n\ge n_1:= {C\over V(t,x,y,z)-\bar V_{n_0}(t,x,y,z)}\vee n_0$. Now assume $n\ge n_1$, and let $\t^n_1>t$ be the first jump time of $Z^n$. Again by
Proposition \ref{smalllarge} (ii) and following similar arguments as in (iii),
we have $|Z^n_{\t^n_1}-z|\ge {1\over n_0}$ on $\{\t^n_1<T\}$. Then,
for any $m$, on $\{\t^n_1<\t_m\}\subset\{\t^n_1<T\}$, using
(\ref{uniconv}) we have
 \beaa
 I^m &\le& I_{\t^n_1} =  V(\t^n_1, X_{\t^n_1}, Y_{\t^n_1}, z)- \bar V_{n_0}(\t^n_1,
 X_{\t^n_1}, Y_{\t^n_1}, z) \\
 &\le& V(\t^n_1, X_{\t^n_1}, Y_{\t^n_1}, z) - V(\t^n_1, X_{\t^n_1}, Y_{\t^n_1}
 -c(Z^n_{\t^n_1}-z), Z^n_{\t^n_1})\\
 &\le& V(\t^n_1, X_{\t^n_1}, Y_{\t^n_1}, z) - V^{n-1}(\t^n_1, X_{\t^n_1}, Y_{\t^n_1}
 -c(Z^n_{\t^n_1}-z), Z^n_{\t^n_1})\\
 &=& V(\t^n_1, X_{\t^n_1}, Y_{\t^n_1}, z) - V^{n}(\t^n_1, X_{\t^n_1}, Y_{\t^n_1},z)\le
 {C\over n}.
 \eeaa
This, together with (\ref{Is}), implies that
\be
\label{tnm}
\lim_{n\to\infty} P(\t^n_1<\t_m) = 0.
\ee

Next, recall from the proof of Theorem \ref{Znexist} that $\t^n_1$ is a
solution to an optimal stopping problem, and thus (cf.
e.g., \cite{ElkarouiKPPQ97}), $V^n(s,X_s,Y_s,z)$ is a martingale for
$t\le \t^n_1$. Therefore
 \beaa
 &&V^n(t,x,y,z) = E\Big\{ V^{n}(\t^n_1\wedge \t_m, X_{\t^n_1\wedge \t_m}, Y_{\t^n_1\wedge
 \t_m}, z)\Big\}\\
 &=& E\Big\{V^{n}(\t_m, X_{\t_m}, Y_{\t_m}, z)\one_{\{\t_m\le \t^n_1\}} + V^{n-1}
  (\t^n_1, X_{\t^n_1}, Y_{\t^n_1}-c(Z^n_{\t^n_1}-z), Z^n_{\t^n_1})\one_{\{ \t^n_1<\t_m\}}
  \Big\}\\
 &\le& E\Big\{V(\t_m, X_{\t_m}, Y_{\t_m}, z)\one_{\{\t_m\le \t^n_1\}}+V(\t^n_1, X_{\t^n_1},
 Y_{\t^n_1}-c(Z^n_{\t^n_1}-z), Z^n_{\t^n_1})\one_{\{ \t^n_1<\t_m\}}\Big\}\\
 &=& E\Big\{V(\t_m, X_{\t_m}, Y_{\t_m}, z) \Big\}\\
 &&+E\Big\{ \Big[V(\t^n_1, X_{\t^n_1}, Y_{\t^n_1}-c(Z^n_{\t^n_1}-z),
 Z^n_{\t^n_1})-V(\t_m, X_{\t_m}, Y_{\t_m},z)]\one_{\{ \t^n_1<\t_m\}}\Big\}.
 \eeaa
Applying Proposition \ref{Vnconv} we then have
$$
V^n(t,x,y,z)\le E\Big\{V(\t_m, X_{\t_m}, Y_{\t_m}, z) + C[1+\sup_{t\le s\le T} |Y_s|]\one_{\{ \t^n_1<\t_m\}}\Big\}.
$$
Sending $n\to\infty$ and by (\ref{tnm}) we obtain (\ref{Vm}), whence the theorem.
\qed

To construct the optimal strategy, we also need

\begin{lem}
\label{outO}
Assume (H1)-(H4). If $(t,x,y)\notin \cO(z)$, then there exists $\tilde z$ such that
$$
V(t,x,y,z) = V(t,x,y-c(\tilde z-z), \tilde z) ~~\mbox{and}~~ (t,x,y-c(\tilde z-z))\in\cO(\tilde z).
$$
\end{lem}

{\it Proof.} Assume the result is not true. Since $(t,x,y)\notin\cO(z)$, there exists $z_1\neq z$ such that $V(t,x,y,z) = V(t,x,y-c(z_1-z), z_1)$. By our assumption, $(t,x,y-c(z_1-z))\notin\cO(z_1)$. Then there exists $z_2\neq z_1$ such that
$V(t,x,y-c(z_1-z), z_1) = V(t,x,y-c(z_1-z)-c(z_2-z_1), z_2)$. Note that $c(z_1-z)+c(z_2-z_1) \ge c(z_2-z)$. By the optimality of $V$ we must have $c(z_1-z)+c(z_2-z_1) = c(z_2-z)$ and
$$
V(t,x,y,z)= V(t,x,y-c(z_1-z),z) =V(t,x,y-c(z_2-z),z_2).
$$
This also implies that $z_2\neq z$.
By our assumption again, $(t,x,y-c(z_2-z))\notin\cO(z_2)$. Repeating 
this argument yields the different $z_1, z_2,\cds$ such that $c(z_i-z)+ c(z_{i+1}-z_i)= c(z_{i+1}-z)$, $i=1, 2,\cds$, and
$$
V(t,x,y,z) = V(t,x,y-c(z_1-z), z_1) =\cds= v(t,x,y-c(z_{i+1}-z), z_{i+1}).
$$
Note that  since $z_i$'s are all different, there is at most one $z_i$ equal to $0$. Thus, by Lemma \ref{cOcompact}, except for one $i$, we have $|z_{i+1}-z_i|\ge \e_1$. This implies that $c(z_i-z) \ge (i-1)[c(\e_1)\wedge c(-\e_1)]$ for all $i$. This contradicts with the fact that $c(z_i-z)$ is bounded.
\qed

We are now ready to construct the optimal strategy $Z^*$. Let $(t,x,y,z)$
be given and denote $X_s := X_s^{t,x}$.

 First,  set $\t^*_0:= t$; if $(t,x,y)\in\cO(z)$, set $Z^*_{\t^*_0} := z$ and $Y^*_{\t^*_0}:= y$; if $(t,x,y)\notin\cO(z)$, applying Lemma \ref{outO} we may find $Z^*_0$ such that
     $
 V(t,x,y,z) = V(t,x,y - c(Z^*_0-z),Z^*_0)
 $ and $(t,x,y - c(Z^*_0-z))\in\cO(Z^*_0)$. In this case, set $Y^*_{\t^*_0} := y-c(Z^*_0-z)$. So in both cases we have $(\t^*_0, X_{\t^*_0}, Y^*_{\t^*_0})\in\cO(Z^*_{\t^*_0})$.

 Assume we have defined $\t^*_i$ and $(Y^*, Z^*)$ on $[t,\t^*_i]$ such that $(\t^*_i, X_{\t^*_i}, Y^*_{\t^*_i})\in \cO(Z^*_{\t^*_i})$. Denote $Y^i_s:= Y^*_{\t^*_i}+Z^*_{\t^*_i}[X_s-X_{\t^*_i}]$, $s\ge \t^*_i$, and
 define
  $$
 \t^*_{i+1}:=\inf \{s\ge \t^*_i: (s,X_s, Y^i_s)\notin\cO(Z^*_{i})\}\wedge T.
 $$
 By Theorem \ref{tnlimit}, $\t^*_{i+1}$ is a stopping time and $\t^*_{i+1}>\t^*_i$ whenever $\t^*_i<T$. Set $Z^*_s := Z^*_{\t^*_i}$ and $Y^*_s := Y^i_s$ for $s\in [\t^*_i, \t^*_{i+1})$. If $\t^*_{i+1}=T$, then we set $ Z^*_{\t^*_{i+1}}:= 0$ and $Y^*_{\t^*_{i+1}} :=
 Y^i_{\t^*_{i+1}} - c(-Z^*_{\t^*_i})$. If $\t^*_{i+1}<T$, by Theorem \ref{tnlimit} again we know $(\t^*_{i+1}, X_{\t^*_{i+1}}, Y^i_{\t^*_{i+1}})\notin \cO(Z^*_{\t^*_i})$. Applying Lemma \ref{outO} we may find $Z^*_{\t^*_{i+1}}$ such that, by defining $Y^*_{\t^*_{i+1}}:= Y^i_{\t^*_{i+1}}- c(Z^*_{\t^*_{i+1}}-Z^*_{\t^*_i})$,
     $$
  V(\t^*_{i+1},X_{\t^*_{i+1}}, Y^i_{\t^*_{i+1}},Z^*_{\t^*_i}) = V(\t^*_{i+1},X_{\t^*_{i+1}}, Y^*_{\t^*_{i+1}},Z^*_{\t^*_{i+1}}),~~\mbox{and}~~(\t^*_{i+1}, X_{\t^*_{i+1}}, Y^*_{\t^*_{i+1}})\in \cO(Z^*_{\t^*_{i+1}}).
 $$

 Repeat the procedure we obtain $\t^*_i$ for $i=0,1,\cds$ and $(Y^*, Z^*)$.

We should point out that at this point we do not know if the above construction will stop after finitely many times.  We shall prove this and our main result Theorem \ref{thm-main} in  Section 7.

\section{Some Technical Proofs}
\setcounter{equation}{0}

In this section we provide the technical proofs we miss in the previous sections. We note that these results are
instrumental in the construction of the piecewise constant optimal
strategy, and some of these results are of interest in their own
right. As a matter of fact, many of these results can be considered
as the necessary conditions of the optimality.

\subsection{Proofs of (\ref{Vnregz}) and Theorem \ref{regzunif}}

To prove the regularity of the $V^n$'s with
respect to $z$, we first introduce the following notion of ``domination" of
strategies. Assume  $Z^j\in \cZ_t^n(z_j)$, $j=1,2$, where either $z_1>
z_2>0$, or $z_1<z_2<0$.  Denote $\D Z:= Z^1-Z^2$, as usual. We say
that $Z^1$ dominates $Z^2$ if $Z^1$ and $Z^2$ have the same jump
times $\t_i$'s, and
 \be
 \label{Dz}
 \D z=\Delta Z_{\tau_{-1}} \ge \Delta Z_{\tau_0} \ge \ldots \ge \Delta Z_{\tau_n} = 0 ~~ \mbox{or}~~ \D z=\Delta Z_{\tau_{-1}} \le \Delta Z_{\tau_0} \le \ldots \le \Delta Z_{\tau_n} = 0,
  \ee
  and, by denoting $\sgn(0):= 0$ and $\d Z^j_{\t_i}:= Z^j_{\t_i}-Z^j_{\t_{i-1}}$,
  \be
  \label{sign}
  \sgn(\d Z^1_{\tau_i}) = \sgn (\d Z^2_{\tau_i}).
   \ee
\begin{rem}
{\rm We remak that the requirements (\ref{Dz}) and (\ref{sign})
guarantee not only that $Z^1$ and $Z^2$ stay close, but that they
are on the same side of the origin. This is mainly due to the fact
that the cost function $c$ is allowed to behave differently
on the two sides of the origin (i.e., $c(0+)\neq c(0-)$).
\qed}
\end{rem}

 Recall (\ref{rhon}). Note that if $z^1>z^2>0$ and
$Z^1$ dominates $Z^2$, then, denoting $Y^i:= Y^{t,x,y,z_i,Z^i}$, $i=1,2$,
and $X=X^{t,x}$,
by induction one can easily check that
 \bea
 \label{dom}
 \Big|E\{U(Y^1_T)\} - E\{U(Y^2_T)\}\Big|
 &\le& C E\Big\{\Big|\int_t^T \Delta Z_s d X_s\Big|
 +\Big|\sum_{i=0}^n [c(\d Z^1_{\t_i}) - c(\d Z^2_{\t_i})]\Big|\Big\}
 \nonumber\\
 &\le& C
|\D z|+ C E\Big\{\sum_{i=0}^n \rho(|\d Z^1_{\tau_i} - \d Z^2_{\tau_i}|)\Big\} \\
 &=& C|\D z| + CE\Big\{\sum_{i=0}^n \rho(\Delta Z_{\tau_{i-1}} -\Delta Z_{\tau_i} )
 \Big\}\le C|\D z| + C\rho_n(|\D z|),
 \nonumber
 \eea

  {\it Proof of (\ref{Vnregz}).}  By the definitions one can easily check that
  \bea
 \label{VT}
 V^n(T,x,y,z) = V(T,z,y,z) = U(y-c(-z)).
 \eea
   Then the estimate is obvious for $t=T$. So we may assume $t<T$.
  Without loss of generality assume $z_1> z_2>0$.

We first prove the right inequality. In light of the estimate (\ref{dom}), we need only prove the
following claim: For any $Z^1\in\cZ_t^n(z_1)$, there exists
$Z^2\in\cZ_t^n(z_2)$ dominated by $Z^1$.
Indeed,  for any $\e>0$, we can
find $Z^{1,\e}\in \cZ^n_t(z_1)$ such that $E\{U(Y^{t,x,y,z,Z^{1,\e}}_T)\}>
V^n(t,x,y,z_1)-\e$. If the claim is true, then (\ref{dom}) leads to that
  $$
  V^n(t,x,y,z_1)\le C[|\D z|+\rho_n(|\D z|)] +V^n(t,x,y,z_2)+\e.
  $$
Letting $\e\to0$ we obtain the right inequality.

Now let $Z^1=\sum_{i=0}^{n-1} Z^1_{\tau_i} \one_{[\tau_i,
\tau_{i+1})}\in \cZ^n_t(z_1)$ be given. We construct $Z^2\in
\cZ^n_t(z_2)$ as follows. We begin by choosing the same jump times
$\t_i$'s. Define
$$
Z^2_{\t_0}:= \left\{\ba{lll} z_2,\q \mbox{if}~ Z^1_{\t_0}=z_1;\\
Z^1_{\t_0},\q\mbox{if}~ Z^1_{\t_0}>z_1 ~\mbox{or}~ Z^1_{\t_0}<z_2;\\
z_2 - {1\over 2}[(z_1-Z^1_{\t_0})\wedge z_2],\q\mbox{if}~ z_2\le Z^1_{\t_0}<z_1.
\ea\right.
$$ Suppose that we have
defined $Z^2_{\t_i}$ such that either $Z^2_{\t_i}=Z^1_{\t_i}$ or
$0<Z^2_{\t_i}<Z^1_{\t_i}$, we then define $Z^2_{\t_{i+1}}$ in the
following way: if $\t_{i+1}=T$ or $Z^2_{\t_i}=Z^1_{\t_i}$, then
simply set $Z^2_{\t_{i+1}}:= Z^1_{\t_{i+1}}$. Assume $\t_{i+1}<T$
and $0<Z^2_{\t_i}<Z^1_{\t_i}$.  Note that in this case, by (\ref{ti})
we have $Z^1_{\t_{i+1}}\neq Z^1_{\t_i}$. If $Z^1_{\tau_{i+1}} >
Z^1_{\tau_i}$ or $Z^1_{\tau_{i+1}} < Z^2_{\tau_i}$, define
$Z^2_{\tau_{i+1}} := Z^1_{\tau_{i+1}}$. Otherwise, we have
$Z^1_{\tau_i} > Z^1_{\tau_{i+1}} \ge Z^2_{\tau_i}>0$, then define
$Z^2_{\tau_{i+1}} := Z^2_{\tau_i} -{1\over 2}[(Z^1_{\t_i}-Z^1_{\t_{i+1}})\wedge Z^2_{\t_i}]$.
Note that we still have either $Z^2_{\t_{i+1}}=Z^1_{\t_{i+1}}$ or
$0<Z^2_{\t_{i+1}}<Z^1_{\t_{i+1}}$, so we may continue to define
$Z^2$. One can check directly that $Z^2$ constructed in such a way satisfies both
(\ref{Dz}) and (\ref{sign}), hence $Z^1$ dominates $Z^2$.

It remains to prove the left inequality. To this end, let $Z^2= \sum_{i=0}^{n-1} Z^2_{\tau_i} \one_{[\tau_i,
\tau_{i+1})}\in \cZ_t^n(z_2)$ be arbitrarily chosen. We define
$Z^1\in\cZ_t^n(z_1)$ recursively
as follows. First, define
$$
Z^1_{\t_0}:= \left\{\ba{lll} z_1,\q \mbox{if}~ Z^2_{\t_0}=z_2;\\
Z^2_{\t_0},\q\mbox{if}~ Z^2_{\t_0}>z_1 ~\mbox{or}~ Z^2_{\t_0}<z_2;\\
z_1 + [Z^2_{\t_0}- z_2],\q\mbox{if}~ z_2< Z^1_{\t_0}\le z_1.
\ea\right.
$$ Assume we have defined
$Z^1_{\t_i}$ such that either $Z^1_{\t_i}=Z^2_{\t_i}$ or
$0<Z^2_{\t_i}<Z^1_{\t_i}$. If $\t_{i+1}=T$ or
$Z^1_{\t_i}=Z^2_{\t_i}$, define $Z^1_{\t_{i+1}}:=
Z^2_{\t_{i+1}}$. Now assume $\t_{i+1}<T$ and
$0<Z^2_{\t_i}<Z^1_{\t_i}$.  Note that in this case
$Z^2_{\t_{i+1}}\neq Z^2_{\t_i}$. If $Z^2_{\tau_{i+1}} <
Z^2_{\tau_i}$ or $Z^2_{\tau_{i+1}} > Z^1_{\tau_i}$, define
$Z^1_{\tau_{i+1}} := Z^2_{\tau_{i+1}}$. Otherwise, we have
$Z^1_{\tau_i} \ge Z^2_{\tau_{i+1}} >Z^2_{\tau_i}>0$, then define
$Z^1_{\tau_{i+1}} := Z^1_{\tau_i} + [Z^2_{\tau_{i+1}} -
Z^2_{\tau_i}]$. Note that we still have either
$Z^1_{\t_{i+1}}=Z^2_{\t_{i+1}}$ or
$0<Z^2_{\t_{i+1}}<Z^1_{\t_{i+1}}$, so we may continue to define
$Z^1$. One may check that (\ref{sign}) still holds, and for each
$\o$, there exists $k$ such that
\bea
\label{DZk}
\Delta Z_{\tau_0} = \cds=\D Z_{\t_k}=\D z~~\mbox{ and}~~ \D Z_{\t_{k+1}}=\cds= \Delta Z_{\tau_n} = 0.
\eea
Then, similar to (\ref{dom}), we have
\beaa
&& E\{U(Y^{t,x,y,z_2, Z^2})\}-V^n(t,x,y,z_1)\le E\{U(Y^{t,x,y,z_2, Z^2})\}-E\{U(Y^{t,x,y,z_1, Z^1})\}\\
\le&& C|\D z|+CE\Big\{\sum_{i=0}^n \rho(\D Z_{\t_{i-1}}-\D Z_{\t_i})\Big\} = C[|\D z|+\rho(|\D z|)].
\eeaa
Since $Z^2$ is arbitrary, we prove the left inequality in (\ref{Vnregz}).
\qed

 {\it Proof of Theorem \ref{regzunif}.} Without loss of generality, assume $z_1>z_2>0$. We
first recall the left inequality in (\ref{Vnregz}).  So we need only check the other half of the inequality. To this end,
let $Z^1$ be the optimal strategy of $V^n(t,x,y,z_1)$, and  as in
the proof of (\ref{Vnregz}) we define $Z^2\in\cZ^n_t(z_2)$ that is ``dominated"
by $Z^1$. We note that, for $i> N(Z^1)$,
$Z^1_{\t_i}=Z^1_{\t_{i-1}}$, which implies that
$Z^2_{\t_i}=Z^2_{\t_{i-1}}$. Then, following (\ref{dom}) we have
 \beaa
  V^n(t,x,y,z_1)-V^n(t,x,y,z_2)&\le& E\{U(Y^{t,x,y,z_1,Z^1}_T)\} - E\{U(Y^{t,x,y,z_2,Z^2}_T)\} \\
 &\le& C|\D z| + CE\Big\{\sum_{i=0}^{N(Z^1)}\rho(|\d Z^1_{\tau_{i}} -\d Z^2_{\tau_i}|)\Big\}\\
 &\le& C|\D z| + C\rho(|\D z|) E\{N(Z^1)\}\le C[|\D z|+\rho(|\D z|)],
 \eeaa
 where the last inequality is due to Theorems \ref{EN1} and
 \ref{EN2}. This proves (\ref{uniregz}).

 To prove (\ref{uniconv}),we denote, for any $m>n$,
 $Z^m=\sum_{i=1}^m Z^m_{\t_{i-1}}\one_{[\t_{i-1},\t_i)}$ be the optimal strategy of
$V^m(t,x,y,z)$. Define $Z^{n,m}_s := Z^m_s
\one_{\{s<\t_n\}}$. Then $Z^{n,m}\in\cZ^{n+1}_t(z)$, and
  \beaa
 Y^{t,x,y,z,Z^m}_T-Y^{t,x,y,z,Z^{n,m}}_T &=& \Big[\int_{\t_n}^T Z^m_s dX_s + c(-Z^m_{\t_{n-1}}) - \sum_{i=n}^m c(\d Z^m_{\t_i})\Big]\one_{\{\t_n<T\}}\\
 &\le& \Big[\int_{\t_n}^T Z^m_s dX_s + c(-Z^m_{\t_{n-1}})\Big]\one_{\{\t_n<T\}}.
 \eeaa
 Note that $\{\t_n<T\}=\{N(Z^m)>n\}$, it follows that
  \beaa
 && V^m(t,x,y,z)-V^{n+1}(t,x,y,z) \le E\Big\{U(Y^{t,x,y,z,Z^m}_T)-U(Y^{t,x,y,z,Z^{n,m}}_T)\Big\}\\
 && \le CE\Big\{\Big[E_{\t_n}\{|\int_{\t_n}^T Z^m_s dX_s|\} + 1\Big]\one_{\{\t_n<T\}}\Big\}\\
 &&\le CP\{\t_n<T\} = CP\{N(Z^m)>n\}\le {C\over n}E\{N(Z^m)\}\le {C\over n}.
 \eeaa
Sending $m\to\infty$ and applying Proposition~\ref{Vnconv}, we obtain
the result.
 \qed
\subsection{Proof of Proposition \ref{smalllarge}}

We split the proof into several lemmas. To begin with, we fix
$(t_0,x_0,y_0,z_0)$ and $n$, and let $Z^n$ be the optimal strategy
of $V^n(t_0,x_0,y_0,z_0)$.
Recall (\ref{Ai}) and for notational simplicity we suppress the superscript ``$n$" and denote them as $A_i$ and $B_i$.
Throughout this subsection we assume that (H1)--(H4) are all in force.
Keep in mind that our purpose is to show that on the set of small
jumps (the set $A_i$'s) the jump will only happen when it jumps
to $0$.
\begin{lem}
\label{Z<z}
 $P$-a.s. on $A_i$, either $0\vee Z^n_{\t_i}\le Z^n_{\t_{i-1}}$ or $Z^n_{\t_{i-1}}\le Z^n_{\t_i}\wedge 0$.
\end{lem}
{\it Proof.}
Suppose that the lemma is not true. Then we may assume
without loss of generality that $P(D_{i_0})>0$ for some $i_0\ge 0$, where
$D_{i_0}:= \{Z^n_{\t_{i_0}}>Z^n_{i_0-1}\ge 0\}\bigcap A_{i_0}$.
Our goal is to construct some $\tilde Z^n\in\cZ^n_{t_0}(z_0)$ such that
\bea
\label{tildeZn>Zn}
E\{U(Y^{\tilde Z^n}_T)\}- E\{U(Y^{Z^n}_T)\}>0, ~~\mbox{where}~Y^{\tilde Z^n}:= Y^{t_0,x_0,y_0,z_0,\tilde Z^n}, Y^{Z^n}:= Y^{t_0,x_0,y_0,z_0, Z^n}.
\eea
This leads to $E\{U(Y^{\tilde Z^n}_T)\}>V(t_0,x_0,y_0,z_0)$, an obvious contradiction.

We now define $\tilde Z^n$ as follows. First, let $k:= \inf\{i\ge
i_0: Z^n_{\t_i}\le 0\}$. Since $Z^n_{\t_n}=0$, we have $k\le n$. Now,
set
  $$
 \tilde Z^n_{\t_i}:=\left\{\ba{lll} Z^n_{\t_i}, & i< i_0~\mbox{or}~i\ge k;\\
 Z^n_{\t_{i_0}-1}\one_{D_{i_0}}+Z^n_{\t_{i_0}}\one_{D_{i_0}^c}, & i=i_0;\\
 \{[Z^n_{\t_i}-Z^n_{\t_{i_0}}+Z^n_{\t_{i_0-1}}]\vee 0\}\one_{D_{i_0}}+Z^n_{\t_i}\one_{D_{i_0}^c}, & i_0+1\le i<k.
 \ea\right.
 $$
Then $\tilde Z^n\in \cZ^n_{\t_0}(z_0)$. To prove (\ref{tildeZn>Zn}), we
denote $\D Z^n:= \tilde Z^n-Z^n$. Then,
$$
\D Y^n_T:= Y^{\tilde Z^n}_T-Y^{Z^n}_T = \int_{\t_{i_0}}^T \D Z^n_s dX_s + \sum_{i=i_0}^n [c(\d Z^n_{\t_i})- c(\d\tilde Z^n_{\t_i})].
$$

By definition of $\tilde Z^n$ it is clear that $\D Y^n_T=0$ on
$D_{i_0}^c$. On $D_{i_0}$, first note that $|\D
Z^n_{\t_i}|\le \d Z^n_{\t_{i_0}}$ for all $i$. Further, for $i>k$, one
has $\d \tilde Z^n_{\t_i}=\d Z^n_{\t_i}$; and for $i\le k$, one can
check  that either $0\le \d\tilde Z^n_{\t_i}\le
\d Z^n_{\t_i}$ or $\d Z^n_{\t_i}\le \d\tilde
Z^n_{\t_i}\le 0$. It then follows from
the monotonicity  assumption in (H3)-(ii) that $c(\d Z^n_{\t_i})\ge c(\d\tilde
Z^n_{\t_i})$. Moreover, note that when $i=i_0$,
$$
c(\d Z^n_{\t_{i_0}})- c(\d\tilde Z^n_{\t_{i_0}}) = c(\d Z^n_{\t_{i_0}})> C_0|\d Z^n_{\t_{i_0}}|,
$$
thanks to Lemma \ref{lem-concave} (i). Thus, on $D_{i_0}$,
$$
\D Y^n_T \ge \int_{\t_{i_0}}^T \D Z^n_s dX_s + c(\d Z^n_{\t_{i_0}})> \int_{\t_{i_0}}^T \D Z^n_s dX_s + C_0|\d Z^n_{\t_{i_0}}|;
$$
and
\beaa
&&
E_{\t_{i_0}}\Big\{|\int_{\t_{i_0}}^T \D Z^n_sdX_s|\Big\}\le
E_{\t_{i_0}}\Big\{\int_{\t_{i_0}}^T |\D Z^n_s b(s,X_s)|ds
+|\int_{\t_{i_0}}^T \D Z^n_s \si(s,X_s)dW_s|\Big\}\nonumber\\
&&\le  E_{\t_{i_0}}\Big\{\int_{\t_{i_0}}^T |\D Z^n_s b(s,X_s)|ds\Big\} + \L E_{\t_{i_0}}\Big\{\int_{\t_{i_0}}^T
 |\D Z^n_s \si(s,X_s)|^2ds\Big\}^{1\over 2}\nonumber\\
&&\le [\|b\|_\infty T +\|\si\|_\infty\sqrt{T}]|\d Z^n_{\t_{i_0}}| =
{\l\over \L} (C_0-1)|\d Z^n_{\t_{i_0}}|.
\eeaa
Therefore, for some appropriately defined $\mathcal F_{T}$-measurable random 
   variable $\xi$, we have
 \beaa
&& E\{U(Y^{\tilde Z^n}_T)-U(Y^{Z^n}_T)\}= E\{U'(\xi)\D Y^n_T\} = E\Big\{U'(\xi)
\D Y^n_T \one_{D_{i_0}}\Big\}\nonumber\\
\ge &&E\Big\{\Big[\l C_0|\d Z^n_{\t_{i_0}}| - \L|\int_{\t_{i_0}}^T \D Z^n_sdX_s|\Big]
 \one_{D_{i_0}}\Big\}\\
 =&& E\Big\{\Big[\l C_0|\d Z^n_{\t_{i_0}}| - \L E_{\t_{i_0}}
 \{|\int_{\t_{i_0}}^T \D Z^n_sdX_s|\}\Big] \one_{D_{i_0}}\Big\}\ge \l
E\{|\d Z^n_{\t_{i_0}}|\one_{D_{i_0}}\}>0.\nonumber
 \eeaa
This proves (\ref{tildeZn>Zn}) and hence the lemma.
\qed

\begin{lem}
\label{smallDz1}
For any $\tilde A_i\subset A_i$, if $P(\tilde A_i)>0$, then
$P(\tilde D_{i+1})>0$, where
 \bea
 \label{PBi}
\tilde D_{i+1}:= \{-1\le {\d Z^n_{\t_{i+1}}\over
\d Z^n_{\t_i}}\le {1\over 2}\}\bigcap\tilde A_i.
 \eea
Consequently, $P$-a.s. in $A_i$, it holds that $|Z^n_{\t_i}|\le
|\d Z^n_{\t_i}|$.
\end{lem}
{\it Proof.} To simplify the presentation we prove the lemma only for $i=1$. The general case can be proved in a line by line analogy. We will prove by contradiction, and without loss of generality, we assume $Z^n_{\t_0}\ge 0$. Then by Lemma \ref{Z<z}, we have $Z^n_{\t_1}< Z^n_{\t_{0}}$ in
$\tilde A_{1}\subset A_{1}$. Suppose that the result is not
true, namely $P(\tilde D_{2})=0$. Then, with possibly an exception of
a null set, one has
$$
\tilde A_1 \subseteq \tilde D_{21}\bigcup \tilde D_{22}:=
(\{\d Z^n_{\t_2}> -\d Z^n_{\t_1}\}\cap A_1)
\bigcup(\{\d Z^n_{\t_2}< {1\over 2}\d Z^n_{\t_1}\}\cap
A_1).
$$
Slightly different from the previous lemma, we now define $\tilde
Z^n_{\t_0}:= Z^n_{\t_0}$; $\tilde Z^n_{\t_1}:= Z^n_{\t_1}\one_{\tilde
A_1^c}+
 z_0\one_{\tilde A_1}$; and $\tilde Z^n_{\t_i}:= Z^n_{\t_i}$, for $i\ge 2$.
 Then $\tilde Z^n\in\cZ^n_{t_0}(z_0)$, and
 $$
 \D Y^n_T = \Big[-\d Z^n_{\t_1}[X_{\t_2}-X_{\t_1}]+ c(\d Z^n_{\t_1})+c(\d Z^n_{\t_2}) - c(Z^n_{\t_2}-Z^n_{\t_0})\Big]\one_{\tilde A_1}.
 $$
Note that, on $\tilde D_{21}$, $Z^n_{\t_2}>Z^n_{\t_0}>Z^n_{\t_1}$. Then
(H3)-(ii) and Lemma \ref{lem-concave} (i) yield that
$$
  c(\d Z^n_{\t_1})+c(\d Z^n_{\t_2}) - c(Z^n_{\t_2}-Z^n_{\t_0})\ge c(\d Z^n_{\t_1})
  \ge C_0|\d Z^n_{\t_1}|.
  $$
On the set $\tilde D_{22}$, however, one has $\d Z^n_{\t_2}<
{1\over 2}\d Z^n_{\t_1}<0$. Then by Lemma \ref{lem-concave} (iii) we have
  $$c(\d Z^n_{\t_1})+c(\d Z^n_{\t_2}) - c(Z^n_{\t_2}-Z^n_{\t_0})\ge
 C_1|\d Z^n_{\t_1}|\ge C_0|\d Z^n_{\t_1}|.
  $$
  So, $P$-a.s. in $\tilde A_1$,
  $$
  \D Y^n_T \ge -\d Z^n_{\t_1}[X_{\t_2}-X_{\t_1}]+
  C_0|\d Z^n_{\t_1}|.
  $$
  Thus, following similar arguments as in Lemma \ref{Z<z}, we have
  \bea
   \label{DY2}
 E\Big\{U(Y^{\tilde Z^n}_T)-U(Y^{Z^n}_T)\Big\} \ge E\Big\{\Big[\l
 C_0|\d Z^n_{\t_1}|
  - \L|\d Z^n_{\t_1}||X_{\t_2}-X_{\t_1}|\Big]\one_{\tilde A_1}\Big\}\ge \l E\Big\{|\d Z^n_{\t_1}|\one_{\tilde A_1}\Big\}>0,
  \eea
a contradiction. Hence $P(\tilde D_2)>0$ must hold.

To prove the last assertion we again assume $i=1$ and $Z^n_{\t_0}\ge 0$,
and that the result is not true. That is, denoting $\h D_1:=
\{|Z^n_{\t_1}|> |\d Z^n_{\t_1}|\}\bigcap A_1$, one has $P(\h
D_1)>0$. Now, denote
$$
\h D_{i+1}:= \{-1\le {\d Z^n_{\t_{i+1}}\over
\d Z^n_{\t_i}}\le {1\over 2}\}\bigcap\h D_i,\q
i=1,\cds,n-1.
$$
We shall prove by induction that that $\h D_i\subset A_i$ and
$Z^n_{\t_{i-1}}\ge Z^n_{\t_i}>{1\over 2}Z^n_{\t_{i-1}}$ on $\h D_i$,
for $i=1,\cds, n$. Indeed,
for $i=1$, by definition $\h D_1\subset A_1$. Moreover, Lemma
\ref{Z<z} tells us that $Z^n_{\t_1}<Z^n_{\t_0}$ on $\h D_1$. If
$Z^n_{\t_1}\le 0$, then obviously $|Z^n_{\t_1}|\le
|\d Z^n_{\t_1}|$. If $Z^n_{\t_1}>0$ in $\h D_1$, then
$Z^n_{\t_1}>-\d Z^n_{\t_1}$ and hence $Z^n_{\t_1}>{1\over 2}Z^n_{\t_0}$ on
$\h D_1$. Namely the claim holds for $i=1$.

Assume now that for all $i\le j$, the claim holds. In particular,
this implies that $Z^n_{\t_j}> {1\over 2^j}Z^n_{\t_0}\ge 0$ on $\h D_j$, we show that the
claim is true for $i=j+1$. Note that on $\h D_{j+1}$, one has
$|\d Z^n_{\t_{j+1}}|\le |\d Z^n_{\t_{j}}|<\e_1$.
Since $Z^n_{\t_j}\neq 0$ on $\hat D_{j+1}\subset \hat D_j$, by (\ref{ti}) 
we know $Z^n_{\t_{j+1}}\neq Z^n_{\t_j}$. Thus $\h D_{j+1}\subset
A_{j+1}$. Moreover, since $\d Z^n_{\t_j}<0$, we have
$\d Z^n_{\t_{j+1}}\geq{1\over 2}\d Z^n_{\t_j}$
on $\h D_{j+1}$. Thus by inductional hypothesis we have
$$
Z^n_{\t_{j+1}}\geq{3\over 2}Z^n_{\t_j} - {1\over 2}Z^n_{\t_{j-1}}
>{1\over 2} Z^n_{\t_j}, \qq \mbox{on $\h D_{j+1}$}.
 $$
That is, the claim is true for $i=j+1$, and hence it is true for all
$i$.

Finally, by applying the same argument
repeatedly we have $P(\h D_n)>0$. But the claim tells us that
$Z^n_{\t_n}> {1\over 2^n}Z^n_{\t_0}\ge 0$ on $\h D_n$. This is impossible
since $Z^n_{\t_n}=0$ must hold almost surely by definition of
$\cZ^n_t(z_0)$. The proof is now complete. \qed

[{\it Proof of Proposition \ref{smalllarge}}] (i) We follow the
arguments in Lemma \ref{smallDz1}. Again for simplicity we assume $i=1$, $Z^n_{\t_0}\ge 0$, and
that the result is not true. Then $P(D_1)>0$, where
$$
D_1:= \Big\{P\{B_2|\cF_{\t_1}\} > {C_0\over C_1}\Big\} \bigcap A_1.
 $$
As before, we define $ \tilde Z^n_{\t_0}:= Z^n_{\t_0}$; $\tilde
Z^n_{\t_1}:= Z^n_{\t_1}\one_{D_1^c}+ Z^n_{\t_0}\one_{D_1}$, and $\tilde Z^n_{\t_i}:= Z^n_{\t_i}$, for $i\ge 2$.
 Then $\tilde Z^n\in\cZ^n_{t_0}(z_0)$, and
 $$
 \D Y^n_T = \Big[-\d Z^n_{\t_1}[X_{\t_2}-X_{\t_1}]+ c(\d Z^n_{\t_1})+c(\d Z^n_{\t_2}) - c(Z^n_{\t_2}-Z^n_{\t_0})\Big]\one_{D_1}.
 $$

 On $D_1\cap B_2^c$, we use (\ref{subadd}) to get
  $
 c(\d Z^n_{\t_1})+c(\d Z^n_{\t_2}) - c(Z^n_{\t_2}-Z^n_{\t_0})\ge 0.
  $
On $D_1\bigcap B_2$, we have $|\d Z^n_{\t_2}|\geq
\e_1\ge |\d Z^n_{\t_1}|.$ Thus Lemma \ref{lem-concave} (iii) tells us that
 $$
 c(\d Z^n_{\t_1})+c(\d Z^n_{\t_2}) - c(Z^n_{\t_2}-Z^n_{\t_0}) \ge C_1|\d Z^n_{\t_1}|.
  $$
Combining above we conclude that
 \bea \label{DY3}
E\{U(Y^{\tilde Z^n}_T)-U(Y^{Z^n}_T)\}
&\ge& E\Big\{\l C_1|\d Z^n_{\t_1}| \one_{D_1\bigcap B_2} - \L|\d Z^n_{\t_1}|E_{\t_1}\{|X_{\t_2}-X_{\t_1}|\}\one_{D_1}\Big\}\nonumber\\
&=& E\Big\{\l C_1|\d Z^n_{\t_1}| E_{\t_1}\{\one_{B_2}\}\one_{D_1} -
\l (C_0-1)|\d Z^n_{\t_1}|\one_{D_1}\Big\}\\
 &\ge& \l E\Big\{\Big[ C_1|\d Z^n_{\t_1}|{C_0\over C_1}-(C_0-1)|\d Z^n_{\t_1}|\Big]
 \one_{D_1}\Big\}=\l E\Big\{|\d Z^n_{\t_1}|\one_{D_1}\Big\}> 0.\nonumber
 \eea
 This is a contradiction and thus proves the part (i).

We shall prove part (ii) by backward induction on $i$. Since
$Z^n_{\t_n}=0$, the result is true for $i=n$. Without loss of
generality we assume it is true for $i=2$ and will prove it for
$i=1$. Assume $Z^n_{\t_0}\ge 0$. If it is not true for $i=1$,
 then $P(\tilde D_1)>0$ where
 $$
 \tilde D_1:= \{Z^n_{\t_1}\neq 0\}\bigcap A_1.
  $$
 We now define $\tilde Z^n_{\t_0}:= Z^n_{\t_0}$; $\tilde Z^n_{\t_1}:= Z^n_{\t_1}
 \one_{\tilde D_1^c}$; and $\tilde Z^n_{\t_i}:= Z^n_{\t_i}$, for $i\ge 2$.
Then $\tilde Z^n\in\cZ^n_{t_0}(z_0)$, and
 $$
 \D Y^n_T = \Big[-Z^n_{\t_1}[X_{\t_2}-X_{\t_1}]+ c(\d Z^n_{\t_1})+c(\d Z^n_{\t_2})-c(-Z^n_{\t_0})-c(Z^n_{\t_2})\Big]1_{\tilde D_1}.
$$
We claim that:
\bea
\label{concave-claim}
\left.\ba{lll}
\dis c(-Z^n_{\t_0}) - c(\d Z^n_{\t_1}) \le {1\over \a_1}[c(-Z^n_{\t_1}) + c(\d Z^n_{\t_1})-c(-Z^n_{\t_0})];\\
\dis c(Z^n_{\t_2}) - c(\d Z^n_{\t_2}) \le {1\over \b_1}[c(-Z^n_{\t_1}) + c(\d Z^n_{\t_1})-c(-Z^n_{\t_0})] + L_0 |Z^n_{\t_1}|;
\ea\right.
~~\mbox{on}~\tilde D_1\bigcap B_2.
\eea
Indeed, without loss of generality, we assume $Z^n_{\t_0}>0$. Then, by Lemmas \ref{Z<z} and \ref{smallDz1}, we have $0\le Z^n_{\t_1} \le -\d Z^n_{\t_1} <\e_1$ on $\tilde D_1 \subset A_1$. Thus the first inequality of (\ref{concave-claim}) follows from Lemma \ref{lem-concave} (ii).  To show the second inequality, note that $|\d Z^n_{\t_2}|\ge \e_1$ on $\tilde D_1\cap B_2$. If $\d Z^n_{\t_2} \le -\e_1$, then $\d Z^n_{\t_2} \le Z^n_{\t_2} <0$, and thus $c(Z^n_{\t_2}) - c(\d Z^n_{\t_2}) \le 0$. If $\e_1\le \d Z^n_{\t_2} \le \e_0$, note that $c(Z^n_{\t_2}) - c(\d Z^n_{\t_2})= c(\d Z^n_{\t_2} + Z^n_{\t_1}) - c( \d Z^n_{\t_2})$ is decreasing in $\d Z^n_{\t_2}$. Then
 \beaa
 c(Z^n_{\t_2}) - c(\d Z^n_{\t_2}) \le c(-\d Z^n_{\t_1} + Z^n_{\t_1}) - c(-\d Z^n_{\t_1}) \le {1\over \b_1}[c(-Z^n_{\t_1}) + c(\d Z^n_{\t_1})-c(-Z^n_{\t_0})],
 \eeaa
 thanks  again to Lemma \ref{lem-concave} (ii). Finally, if  $\d Z^n_{\t_2}\ge \e_0$, then
 $ c(Z^n_{\t_2}) - c(\d Z^n_{\t_2})  \le L_0|Z^n_{\t_1}|$.
 This completes the proof of Claim (\ref{concave-claim}).

Note that,  by inductional hypothesis we have
$Z^n_{\t_2}=0$ on $\tilde D_1\bigcap B_2^c$. Then, for some appropriately defined $\mathcal F_{T}$-measurable random 
   variable $\xi$, by (\ref{subadd}) we have,
 \beaa
&&E\{U(Y^{\tilde Z^n}_T)-U(Y^{Z^n}_T)\} = E\{U'(\xi)[Y^{\tilde Z^n}_T-Y^{Z^n}_T]\}\\
&=& E\Big\{U'(\xi)\Big[-Z^n_{\t_1}[X_{\t_2}-X_{\t_1}]\one_{\tilde D_1} + [c(-Z^n_{\t_1}) + c(\d Z^n_{\t_1})-c(-Z^n_{\t_0})]\one_{\tilde D_1B_2^c}\\
&&\qq + [c(\d Z^n_{\t_1})-c(-Z^n_{\t_0})+c(\d Z^n_{\t_2})-c(Z^n_{\t_2})] \one_{\tilde D_1 B_2}\Big]\Big\}\\
&\ge& E\Big\{U'(\xi)\Big[-Z^n_{\t_1}[X_{\t_2}-X_{\t_1}]\one_{\tilde D_1} - L_0 |Z^n_{\t_1}|  \one_{\tilde D_1 B_2}\\
&&\qq + [c(-Z^n_{\t_1}) + c(\d Z^n_{\t_1})-c(-Z^n_{\t_0})][\one_{\tilde D_1B_2^c}-({1\over\a_1}+{1\over \b_1})
 \one_{\tilde D_1 B_2}]\Big]\Big\}\\
 &\ge& E\Big\{-\L|Z^n_{\t_1}|[|X_{\t_2}-X_{\t_1}|+L_0]\one_{\tilde D_1}\\
 &&+ [c(-Z^n_{\t_1}) + c(\d Z^n_{\t_1})-c(-Z^n_{\t_0})][\l \one_{\tilde D_1B_2^c}
 -\L({1\over\a_1}+{1\over \b_1})\one_{\tilde D_1B_2}]\Big\}\\
 &=& E\Big\{\Big[-\L|Z^n_{\t_1}|[E_{\t_1}\{|X_{\t_2}-X_{\t_1}|\} + L_0]\\
 && + [c(-Z^n_{\t_1}) + c(\d Z^n_{\t_1})-c(-Z^n_{\t_0})][\l E_{\t_1}\{\one_{B_2^c}\}
 -\L ({1\over\a_1}+{1\over \b_1})E_{\t_1}\{\one_{B_2}\}]\Big]\one_{\tilde D_1}\Big\}.
 \eeaa
 One can easily check that
 \beaa
\L \Big[E_{\t_1}\{|X_{\t_2}-X_{\t_1}|\} + L_0\Big] \le \l[C_0-1].
 \eeaa
 Moreover,  by part (i)
we know that $P$-a.s. on $\tilde D_1\subset A_1$, $
P\{B_2|\cF_{\t_1}\}\le {C_0\over C_1}$ and thus $P\{B_2^c|\cF_{\t_1}\}\ge
1-{C_0\over C_1}$. Then
\beaa
\l E_{\t_1}\{\one_{B_2^c}\} -{\L\over\a_1}E_{\t_1}\{\one_{B_2}\}\ge \l\Big[1-{C_0\over C_1} - \L({1\over\a_1}+{1\over \b_1}) {C_0\over C_1}\Big] = {\l C_0\over C_1}.
\eeaa
Note that, on $\tilde D_1 \subset A_1$, by Lemmas \ref{Z<z} and \ref{smallDz1}, we have $0\le Z^n_{\t_1} \le \d Z^n_{\t_1}<\e_1$ or $0>\d Z^n_{\t_1}\ge Z^n_{\t_1}>-\e_1$.  Then it follows from Lemma \ref{lem-concave} (ii) that
 \bea
 \label{DY4}
&&E\{U(Y^{\tilde Z^n}_T)-U(Y^{Z^n_T})\}\\
&\ge&  E\Big\{\Big[-\l(C_0-1)|Z^n_{\t_1}|+ {\l C_0\over C_1}C_1|Z^n_{\t_1}|
\Big]\one_{\tilde D_1}\Big\}= \l
E\Big\{|Z^n_{\t_1}|\one_{\tilde D_1}\Big\}>0,\nonumber
\eea
a contradiction.
 \qed

\bs

\subsection{Proofs of Lemma \ref{cOcompact} and Theorem \ref{thm-main}}

[{\it Proof of Lemma \ref{cOcompact}.}] We follow the proof of
Proposition \ref{smalllarge}.
 For each $n$, let $Z^{n}\in \cZ^{n}_t(\tilde z)$ be the optimal portfolio of
$V^{n}(t,x,y-c(\tilde z-z),\tilde z)$. We first prove several claims
by contradiction. In each case, we show that if the claim is not
true, then we can construct some $\tilde Z^n\in\cZ^{n+1}_t(z)$ such
that, by denoting $Y^{\tilde Z^n}:= Y^{t,x,y,z,\tilde Z^n}, Y^{ Z^n}:= Y^{t,x,y-c(\tilde z-z),\tilde z,Z^n}$,
 \bea
 \label{DV}
 E\Big\{U(Y^{\tilde Z^n}_T)-U(Y^{Z^n}_T)\Big\} \ge c(z,\tilde z)>0
 \eea
 where $c(z,\tilde z)$ is some constant independent of $n$. This implies that
 $$V(t,x,y,z) - V_n(t,x,y-c(\tilde z-z),\tilde z)\ge c(z,\tilde z)>0.
 $$
 Sending $n\to \infty$ and applying Proposition \ref{Vnconv}, we obtain the contradiction.

Without loss of generality we assume $z\ge 0$. The key observation
is that we may also view $Y^{Z^n}_T$ as the wealth of the portfolio $Z^{n}$ starting from
$(t,x,y,z)$, with two initial jumps first from $z$ to $\tilde z$ and then from $\tilde z$ to $Z^n_{\t_0}$.


{\bf Claim 1.} $\tilde z< z$. Indeed, if $\tilde z>z$, for fixed $n$, let
$k:= \inf\{i\ge 0: Z^n_{\t_i}\le 0\}$, and define $\tilde
Z^n_{\t_i}:= [Z^n_{\t_i}-\tilde z+z]\vee 0$ for $i<k$, and
$\tilde Z^n_{\t_i}:= Z^n_{\t_i}$ for $i\ge k$. Then $\tilde
Z^n\in \cZ^{n+1}_t(z)$. Following exactly the same arguments as in the proof of Lemma \ref{Z<z},
we prove (\ref{DV}) with $c(z,\tilde z) = \l(\tilde z-z)>0$ and thus obtain a contradiction.

\ms

{\bf Claim 2.} $-1\le {Z^n_{\t_0}-\tilde z\over \tilde z-z}\le {1\over 2}$, and if $Z^n_{\t_0}=\tilde z$ then
$P(-1\le {Z^n_{\t_1}-\tilde z\over \tilde z-z}\le {1\over 2}) >0$.
Indeed, assume the result is not true. Define $\tilde Z^n_{\t_0}=z$, $\tilde
Z^n_{\t_i}=Z^n_{\t_i}$, for all $i\geq 1$. Then $\tilde Z^n\in
\cZ_t^{n+1}(z)$, and similar to (\ref{DY2}) we prove (\ref{DV}) with $c(z,\tilde z) = \l(z-\tilde z)>0$.


\ms {\bf Claim 3.} $|\tilde z|\le z-\tilde z$. Indeed, if $Z^n_{\t_0}\neq \tilde z$, then by Claim 2 we have $0<|Z^n_{\t_0}-\tilde z| \le z-\tilde z<\e_1$. Applying Proposition
\ref{smalllarge} we get $Z^n_{\t_0}=0$ and thus proving the claim. If $Z^n_{\t_0}= \tilde z$, then Claim 2 leads to $P(-1\le {Z^n_{\t_1}-\tilde z\over \tilde z-z}\le {1\over 2}) >0$. On $\{-1\le {Z^n_{\t_1}-\tilde z\over \tilde
z-z}\le {1\over 2}\}$, we have $|Z^n_{\t_1}-\tilde z|\le |\tilde
z-z|<\e_1$. If $Z^n_{\t_1}=\tilde z$, by (\ref{ti}) we get $\t_1=T$ and thus $\tilde z=0$. If
$Z^n_{\t_1}\neq\tilde z$,  by Proposition
\ref{smalllarge} again we get $Z^n_{\t_1}=0$. Then $-1\le {-\tilde
z\over \tilde z-z}\le {1\over 2}$ and thus the claim holds.

\ms {\bf Claim 4.} If $Z^n_{\t_0}= \tilde z$, then $ P(|Z^n_{\t_1}-\tilde z|\ge \e_1)\le {C_0\over
C_1}$.  Indeed, if $P(|Z^n_{\t_1}-\tilde z|\ge \e_1)> {C_0\over
C_1}$, then we define $ \tilde Z^n_{\t_0}:= z$, and $\tilde
Z^n_{\t_i}:= Z^n_{\t_i}$, for $i\geq 1$. Similar to (\ref{DY3}) we prove (\ref{DV}) with $c(z,\tilde z) = \l(z-\tilde z)>0$.

We now prove the lemma. Define $\tilde Z^n_{\t_0}:= 0$ and $\tilde Z^n_{\t_i}:= Z^n_{\t_i}$ for $i\ge 1$.
Then
\beaa
\D Y_T &:=& Y^{\tilde Z^n}_T - Y^{Z^n}_T \\
&=& c(\tilde z-z)+c(Z^n_{\t_0}-\tilde z) + C(Z^n_{\t_1}-Z^n_{\t_0})-c(-z)-c(Z^n_{\t_1})-Z^n_{\t_0}[X_{\t_1}-x].
\eeaa
If $Z^n_{\t_0}\neq \tilde z$, by the proof of Claim 3, we have $Z^n_{\t_0}=0$. Then
$$
\D Y_T = c(\tilde z-z)+c(-\tilde z)-c(-z)\ge C_1|\tilde z|,
$$
thanks to Lemma \ref{lem-concave} (iii) and Claims 1 and 3. If $Z^n_{\t_0}=\tilde z$, then
$$
\D Y_T = c(\tilde z-z)+c(Z^n_{\t_1}-\tilde z)
-c(-z)- c(Z^n_{\t_1})-\tilde z [X_{\t_1}-x].
$$
Similar to (\ref{DY4}) we can prove
$$
V(t,x,y,z)-V^n(t,x,y-c(\tilde z-z),\tilde z) \ge E\Big\{U(Y^{\tilde Z^n}_T)
-U( Y^{Z^n}_T)\Big\} \ge \l |\tilde z|.
$$
Send $n\to\infty$ and noting that $V(t,x,y,z)=V(t,x,y-c(\tilde
z-z),\tilde z)$, we must have
$\tilde z=0$.
\qed

\bs

\no[{\it Proof of Theorem \ref{thm-main}.}] (i) If $c(z)\ge c_0>0$ for
all $z\neq0$, then following the arguments in Theorem \ref{EN1} one
can easily prove that $E\{N(Z^*)\} \le {C\over c_0}$.

Now assume (H4) holds. Following the proof of Theorem
\ref{cOcompact} it is readily seen that, for any $i\ge 1$ and $P$-a.s. on $\{0<|Z^*_{\t^*_i}-Z^*_{\t^*_{i-1}}|<\e_1\}$, it holds
that
$$
Z^*_{\t^*_i}=0 \q \mbox{and}\q P(|Z^*_{\t^*_{i+1}}-Z^*_{\t^*_i}|\ge \e_1|\cF_{\t^*_i})
\le {C_0\over
C_1}.
$$
Then following the proof of Theorem \ref{smalljump} we get
$$
P\Big(\sum_{i=0}^n \one_{\{0<|Z^*_{\t^*_i}-Z^*_{\t^*_{i-1}}|<\e_1\}}
\ge m\Big) \le {1\over 2^{m}},\q\forall n\ge m.
$$
Similar to Theorem \ref{EN2} one can then prove that
$E\Big\{\sum_{i=0}^\infty
\one_{\{|Z^*_{\t^*_i}-Z^*_{\t^*_{i-1}}|>0\}}\Big\} <\infty$. This
implies that $P(\t^*_i<T,\forall i)=0$ and $E(N(Z^*))<\infty$.

(ii) Applying Lemma \ref{tnlimit} repeatedly we have
$$
V(t,x,y,z) = E\Big\{V(\t^*_n, X_{\t^*_n}, Y^*_{\t^*_n}, Z^*_{\t^*_n})\Big\},\q \forall n.
$$
Now by (i), we conclude that $\t^*_n=T$ and $Z^*_{\t^*_n}=0$ for $n$
large enough. Sending $n\to\infty$ we obtain that $ V(t,x,y,z) =
E\{U(Y^*_T)\}$. This means that $Z^*$ is an optimal portfolio for
$V(t,x,y,z)$. \qed


\end{document}